\documentclass[11pt]{article}%
\usepackage{amssymb}
\usepackage[all]{xy}
\usepackage[latin1]{inputenc}
\usepackage{amsfonts}
\usepackage{amsmath}
\usepackage{amssymb}
\usepackage{url}
\usepackage{hyperref}
\usepackage{graphicx}%
\providecommand{\U}[1]{\protect\rule{.1in}{.1in}}
\numberwithin{equation}{section}
\providecommand{\U}[1]{\protect\rule{.1in}{.1in}}
\providecommand{\U}[1]{\protect\rule{.1in}{.1in}}
\newtheorem{theo}{Theorem}[section]
\newtheorem{prop}[theo]{Proposition}
\newtheorem{lem}[theo]{Lemma}

\begin{document}

\title{Particle approximation of the intensity measures of a spatial branching point
process arising in multi-target tracking}
\author{Francois Caron\thanks{Centre INRIA Bordeaux Sud-Ouest \& Institut de
Mathématiques de Bordeaux , Université Bordeaux, 351 cours de la Libération
33405 Talence cedex, France, Francois.Caron@inria.fr} , Pierre Del
Moral\thanks{Centre INRIA Bordeaux Sud-Ouest \& Institut de Mathématiques de
Bordeaux , Université Bordeaux, 351 cours de la Libération 33405 Talence
cedex, France, Pierre.Del-Moral@inria.fr} , Arnaud Doucet\thanks{Department of
Statistics \& Department of Computer Science, University of British Columbia,
333-6356 Agricultural Road, Vancouver, BC, V6T 1Z2, Canada,
Arnaud@stat.ubc.ca}, Michele Pace\thanks{Centre INRIA Bordeaux Sud-Ouest \&
Institut de Mathématiques de Bordeaux , Université Bordeaux, 351 cours de la
Libération 33405 Talence cedex, France, Michele.Pace@inria.fr} }
\maketitle

\begin{abstract}
The aim of this paper is two-fold. First we analyze the sequence of intensity
measures of a spatial branching point process arising in a multiple target
tracking context. We study its stability properties, characterize its long
time behavior and provide a series of weak Lipschitz type functional
contraction inequalities. Second we design and analyze an original particle
scheme to approximate numerically these intensity measures. Under appropriate
regularity conditions, we obtain uniform and non asymptotic estimates and a
functional central limit theorem. To the best of our knowledge, these are the
first sharp theoretical results available for this class of spatial branching
point processes.


\end{abstract}

\section{Introduction}

Multi-target tracking problems deal with tracking several targets
simultaneously given noisy sensor measurements. Over recent years, point
processes approaches to address these problems have become very popular. The
use of point processes in a multiple-target tracking context was first
proposed in S. Mori et al.~\cite{mori} as early as in 1986. Using a random
sets formalism, a formalism essentially equivalent to the point process
formalism \cite{mori2002}, R. Mahler and his co-authors proposed in two books
~\cite{maler4,maler5} a systematic treatment of multi-sensor multi-target
filtering problems. However, as mentioned in \cite{mori2002},
\textquotedblleft... although the random sets formalism (or the point process
formalism) for multitarget tracking has provided a unified view on the subject
of multiple target tracking, it has failed to produce any significant
practical tracking algorithms...\textquotedblright.

This situation has recently changed following the introduction of the
PHD\ (probability hypothesis density) filter by R. Mahler \cite{maler0,maler}.
The PHD\ filter is a powerful multi-target tracking algorithm which is
essentially a Poisson type approximation to the optimal multi-target filter
\cite{maler0,maler,singh2009}. It has found numerous applications since its
introduction. The PHD\ filter cannot be computed analytically but it can be
approximated by a mixture of Gaussians for linear Gaussian target models
\cite{vo2006} and by non-standard particle methods for nonlinear non-Gaussian
target models \cite{arnaud2,arnaud3}.

Despite their increasing popularity, the theoretical performance of these
multi-target particle methods remain poorly understood. Indeed their
mathematical structure is significantly different from standard particle
filters so the detailed theoretical\ results for particle filters provided in
\cite{fk} are not applicable. Some convergence results have already been
established in \cite{arnaud2,arnaud3} but remain quite limited. Reference
\cite{arnaud2} presents a basic convergence result for the PHD\ filter but
does not establish any rate of convergence. In \cite{arnaud3} the authors
provide some quantitative bounds and a central limit theorem. However these
quantitative bounds are not sharp and no stability result is provided.

The aim of this work is to initiate a thorough theoretical study of these
non-standard particle methods by first characterizing the stability properties
of the \textquotedblleft signal\textquotedblright\ process and establishing
uniform w.r.t. the time index convergence results for its particle
approximation. This \textquotedblleft signal\textquotedblright\ process is a
spatial branching point process whose intensity measure always satisfies a
closed recursive equation in the space of bounded positive measures. We will
not consider any observation process in this article. The analysis of the
particle approximations of PHD filters is presented in \cite{cdpv}. It builds
heavily upon the present work but is even more complex as it additionally
involves at each time step a nonlinear update of the intensity measure.

The rest of this paper is organized as follows:

In section~\ref{secintro}, we present a spatial branching point process which
is general enough to model a wide variety of multiple target problems. We
establish the linear evolution equation associated to the intensity measures
of this process and introduce an original particle scheme to approximate them
numerically. Section \ref{sec:statementresults} summarizes the main results of
this paper. In Section \ref{branchingsectsg}, we provide a detailed analysis
of the stability properties and the long time behavior of this sequence of
intensity measures, including the asymptotic behavior of the total mass
process, i.e. the integral of the intensity measure over the state space, and
the convergence to equilibrium of the corresponding sequence of normalized
intensity measures. For time-homogeneous models, we exhibit three different
types of asymptotic behavior. The analysis of these stability properties is
essential in order to guarantee the robustness of the model and to obtain
reliable numerical approximation schemes. Section \ref{meanfieldsec} is
devoted to the theoretical study of the non-standard particle scheme
introduced to approximate the intensity measures. Our main result in this
section is a non-asymptotic convergence for this scheme. Under some
appropriate stability conditions, we additionally obtain uniform estimates
w.r.t. the time parameter.

\section{Spatial branching point process and its particle
approximation\label{secintro}}

\subsection{Spatial branching point process for multi-target tracking}

\label{secinspa} Assume that at a given time $n$ there are $N_{n}$ target
states $(X_{n}^{i})_{1\leq i\leq N_{n}}$ taking values in some measurable
state space $E_{n}$ enlarged with an auxiliary cemetery point $c$. The state
space $E_{n}$ depends on the problem at hand. It may vary with the time
parameter and can include all the characteristics of a target such as its
type, its kinetic parameters as well as its complete path from the origin. As
usual, we extend the measures $\gamma_{n}$ and the bounded measurable
functions $f_{n}$ on $E_{n}$ by setting $\gamma_{n}(c)=0$ and $f_{n}(c)=0$.

Each target has a survival probability $e_{n}(X_{n}^{i})\in\lbrack0,1]$. When
a target dies, it goes to the cemetery point $c$. We also use the convention
$e_{n}(c)=0$ so that a dead target can only stay in the cemetery. Survival
targets give birth to a random strictly positive number of individuals
$h_{n}^{i}(X_{n}^{i})$ where $\left(  h_{n}^{i}(X_{n}^{i})\right)  _{1\leq
i\leq N_{n}}\ $is a collection of independent random variables such that
$\mathbb{E}\left(  h_{n}^{i}(x_{n})\right)  =H_{n}(x_{n})$ for any $x_{n}\in
E_{n}$ where $H_{n}$ is a given collection of bounded functions $H_{n}$. We
have $H_{n}\left(  x_{n}\right)  \geq1$ for any $x_{n}\in E_{n}$ as $h_{n}%
^{i}(x_{n})\geq1$. This branching transition is called spawning in the
multi-target tracking literature. We define $G_{n}=e_{n}H_{n}$.

After this branching transition, the system consists of a random number
$\widehat{N}_{n}$ of individuals $(\widehat{X}_{n}^{i})_{1\leq i\leq
\widehat{N}_{n}}$. Each of them evolves randomly $\widehat{X}_{n}^{i}%
=x_{n}\leadsto X_{n+1}^{i}$ according to a Markov transition $M_{n+1}%
(x_{n},dx_{n+1})$ from $E_{n}$ into $E_{n+1}$. We use the convention
$M_{n+1}(c,c)=1$, so that any dead target remains in the cemetery state.

At the same time, an independent collection of new targets is added to the
current configuration. This additional point process is modeled by a spatial
Poisson process with a prescribed intensity measure $\mu_{n+1}$ on $E_{n+1}$.
It is used to model new targets entering the state space.

At the end of this transition, we obtain $N_{n+1}=\widehat{N}_{n}%
+N_{n+1}^{\prime}$ targets $(X_{n+1}^{i})_{1\leq i\leq N_{n+1}}$, where
$N_{n+1}^{\prime}$ is a Poisson random variable with parameter given by the
total mass $\mu_{n+1}(1)$ of the positive measure $\mu_{n+1}$, and
$(X_{n+1}^{\widehat{N}_{n}+i})_{1\leq i\leq N_{n+1}^{\prime}}$ are independent
and identically distributed random variables with common distribution
$\overline{\mu}_{n+1}=\mu_{n+1}/\mu_{n+1}(1)$ where $\mu_{n+1}(1):=\int%
_{E_{n+1}}\mu_{n+1}(dx)$.\bigskip

\textbf{Example.}\emph{ To illustrate the model, we present here a simple yet
standard example \cite{vo2006} of a target evolving in a two-dimensional
surveillance region $S\subset\mathbb{R}^{2}$. In this case, we set
$E_{n}=E=S\times\mathbb{R}^{2}$. All the targets are assumed to be of the same
type. The state of a target $X_{n}=\left[  p_{n}^{x},p_{n}^{y},v_{n}^{x}%
,v_{n}^{y}\right]  ^{T}$ consists of its position $\left(  p_{n}^{x},p_{n}%
^{y}\right)  \in S$ and velocity $\left(  v_{n}^{x},v_{n}^{y}\right)
\in\mathbb{R}^{2}$ and is assumed to evolve according to a linear Gaussian
model
\begin{equation}
X_{n}=AX_{n-1}+V_{n}\label{eq:lineargaussian}%
\end{equation}
where $A$ is a known transition matrix and $V_{n}\sim\mathcal{N}\left(
0,\Sigma\right)  $ is a sequence of i.i.d zero-mean normal random variables of
covariance $\Sigma$; i.e. $M_{n}(x_{n-1},dx_{n})=M(x_{n-1},x_{n})dx_{n}$ with
\[
M(x_{n-1},x_{n})=\left\vert 2\pi\Sigma\right\vert ^{-1/2}\exp\left(  -\frac
{1}{2}\left(  x_{n}-Ax_{n-1}\right)  ^{\text{T}}\Sigma^{-1}\left(
x_{n}-Ax_{n-1}\right)  \right)  .
\]
In the example, we assume that $\mu_{n}\left(  x\right)  =\mu\left(  x\right)
$, $e_{n}(x)=s>0$ and $h_{n}(x_{n})=h\in\left\{  1,2\right\}  $ with
$\mathbb{P}\left(  h=1\right)  =1-\mathbb{P}\left(  h=2\right)  =\alpha$.
Hence for this model, each target $X_{n-1}$ survives at time $n-1$ with a
probability $s$. Each survival target has one offspring with probability
$\alpha$ which evolves according to (\ref{eq:lineargaussian}) or two offspring
with probability $1-\alpha$ which, conditional upon $X_{n-1}$, independently
evolve according to (\ref{eq:lineargaussian}). Additionally, a random number
of targets distributed according to a Poisson distribution of parameter
$\mu\left(  1\right)  $ appear. These targets are independent and distributed
in $E$ as $\overline{\mu}=\mu/\mu\left(  1\right)  $.}

\subsection{Sequence of intensity distributions}

\label{intensityflows}

At every time $n$, the intensity measure of the point process $\mathcal{X}%
_{n}:=\sum_{i=1}^{N_{n}}\delta_{X_{n}^{i}}$ associated to the targets is given
for any bounded measurable function $f$ on $E_{n}\cup\{c\}$ by the following
formula:
\[
\gamma_{n}(f):=\mathbb{E}\left(  \mathcal{X}_{n}(f)\right)  \quad\mbox{\rm
with}\quad\mathcal{X}_{n}(f):=\int~f(x)~\mathcal{X}_{n}(dx)
\]
To simplify the presentation, we suppose that the initial configuration of the
targets is a spatial Poisson process with intensity measure $\mu_{0}$ on the
state space $E_{0}$.

Given the construction defined in section~\ref{secinspa}, it follows
straightforwardly that the intensity measures $\gamma_{n}$ on $E_{n}$ satisfy
the following recursive equation.

\begin{lem}
For any $n\geq0,$ we have%
\begin{equation}
\gamma_{n+1}(dx^{\prime})=\int~\gamma_{n}(dx)~Q_{n+1}(x,dx^{\prime})+\mu
_{n+1}(dx^{\prime}) \label{intensity}%
\end{equation}
with the initial condition $\gamma_{0}=\mu_{0}$ where $\mu_{n+1}$ is the
intensity measure of the spatial point process associated to the birth of new
targets at time $n+1$ while the integral operator $Q_{n+1}$ from $E_{n}$ into
$E_{n+1}$ is defined by
\begin{equation}
Q_{n+1}(x_{n},dx_{n+1}):=G_{n}(x_{n})~M_{n+1}(x_{n},dx_{n+1}).
\label{fkmodels}%
\end{equation}

\end{lem}

\noindent\mbox{\bf Proof:}\newline For any bounded measurable function $f$ on
$E_{n+1}\cup\{c\}$, we have
\[
\gamma_{n+1}\left(  f\right)  =\mathbb{E}\left(  \sum_{i=1}^{\widehat{N}_{n}%
}f\left(  X_{n+1}^{i}\right)  \right)  +\mathbb{E}\left(  \sum_{i=\widehat{N}%
_{n}}^{\widehat{N}_{n}+N_{n+1}^{\prime}}f\left(  X_{n+1}^{i}\right)  \right)
\]
where, thanks to the Poisson assumption, we have
\[
\mathbb{E}\left(  \sum_{i=\widehat{N}_{n}}^{\widehat{N}_{n}+N_{n+1}^{\prime}%
}f\left(  X_{n+1}^{i}\right)  \right)  =\mu_{n+1}\left(  1\right)
\overline{\mu}_{n+1}\left(  f\right)  =\mu_{n+1}(f)
\]
and
\begin{align*}
\mathbb{E}\left(  \sum_{i=1}^{\widehat{N}_{n}}f\left(  X_{n+1}^{i}\right)
\right)   &  =\mathbb{E}\left(  \mathbb{E}\left(  \left.  \sum_{i=1}%
^{\widehat{N}_{n}}f\left(  X_{n+1}^{i}\right)  \right\vert \mathcal{F}%
_{n}\right)  \right) \\
&  =\mathbb{E}\left(  \mathbb{E}\left(  \left.  \sum_{i=1}^{\widehat{N}_{n}%
}M_{n+1}\left(  f\right)  \left(  \widehat{X}_{n}^{i}\right)  \right\vert
\mathcal{G}_{n}\right)  \right) \\
&  =\mathbb{E}\left(  \sum_{i=1}^{N_{n}}e_{n}\left(  X_{n}^{i}\right)
H_{n}(X_{n}^{i})M_{n+1}\left(  f\right)  \left(  X_{n}^{i}\right)  \right) \\
&  =\gamma_{n}\left(  e_{n}H_{n}M_{n+1}\left(  f\right)  \right)
\end{align*}
where $\mathcal{F}_{n}$ denotes the $\sigma$-field generated by $(\widehat{X}%
_{n}^{i})_{1\leq i\leq\widehat{N}_{n}}$ and $\mathcal{G}_{n}$ the $\sigma
$-field generated by $\left(  X_{n}^{i}\right)  _{1\leq i\leq N_{n}}$.
\hfill\hbox{\vrule height 5pt width 5pt depth 0pt}\medskip\newline

These intensity measures typically do not admit any closed-form expression. A
natural way to approximate them numerically is to use a particle
interpretation of the associated sequence of probability distributions given
by
\[
\eta_{n}(dx):=\gamma_{n}(dx_{n})/\gamma_{n}(1)\quad\mbox{\rm with}\quad
\gamma_{n}(1):=\int_{E_{n}}\gamma_{n}(dx)
\]
To avoid unnecessary technical details, we further assume that the potential
functions $G_{n}$ are chosen so that for any $x\in E_{n}$
\begin{equation}
0<g_{n,-}\leq G_{n}(x)\leq g_{n,+}<\infty\label{eq:assumptionpotential}%
\end{equation}
for any time parameter $n\geq0$. Note that this assumption is satisfied in
most realistic multi-target scenarios such as the example discussed at the end
of section \ref{secinspa}. Indeed the condition $g_{n,-}\leq G_{n}(x)$
essentially states that there exists $e_{n,-}>0$ such that $e_{n}\left(
x\right)  \geq e_{n,-}$ for any $x\in E_{n}$ as $H_{n}\left(  x\right)  \geq
1$. The condition $G_{n}(x)\leq g_{n,+}$ states that there exists
$H_{n,+}<\infty$ such that $H_{n}\left(  x\right)  \leq H_{n,+}$ for any $x\in
E_{n}$ as $e_{n}\left(  x\right)  \leq1$. In the unlikely scenario where
(\ref{eq:assumptionpotential}) is not satisfied then the forthcoming analysis
can be extended to more general models using the techniques developed in
section 4.4 in~\cite{fk}; see also~\cite{cerou2008}. We denote by
$\mathcal{P}(E_{n})$ the set of probability measures on the state space
$E_{n}$.

To describe these particle approximations, it is important to observe that the
pair process $(\gamma_{n}(1),\eta_{n})\in(\mathbb{R}_{+}\times\mathcal{P}%
(E_{n}))$ satisfies an evolution equation of the following form
\begin{equation}
(\gamma_{n}(1),\eta_{n})=\Gamma_{n}(\gamma_{n-1}(1),\eta_{n-1})
\label{flotPhi}%
\end{equation}
We let $\Gamma_{n}^{1}$ and $\Gamma_{n}^{2}$ be the first and the second
component mappings from $(\mathbb{R}_{+}\times\mathcal{P}(E_{n}))$ into
$\mathbb{R}_{+}$, and from $(\mathbb{R}_{+}\times\mathcal{P}(E_{n}))$ into
$\mathcal{P}(E_{n})$. The mean field particle approximation associated with
the equation (\ref{flotPhi}) relies on the fact that it is possible to rewrite
the mapping $\Gamma_{n+1}^{2}$ in the following form
\begin{equation}
\Gamma_{n+1}^{2}(\gamma_{n}(1),\eta_{n})=\eta_{n}K_{n+1,(\gamma_{n}%
(1),\eta_{n})} \label{flotPhimc}%
\end{equation}
where $K_{n+1,(m,\eta)}$ is a Markov kernel indexed by the time parameter $n
$, a mass parameter $m\in\mathbb{R}_{+}$ and a probability measure $\eta$ on
the space $E_{n}$. In the literature on mean field particle systems,
$K_{n,(m,\eta)}$ is called a McKean transition. The choice of \ such Markov
transitions $K_{n,(m,\eta)}$ is not unique and will be discussed in
section~\ref{nonlinsecc}.

Before concluding this section, we note that
\begin{equation}
\gamma_{n+1}(dx^{\prime})=(\gamma_{n}Q_{n+1})(dx^{\prime}):=\int\gamma
_{n}(dx)~Q_{n+1}(x,dx^{\prime}) \label{CaseNull}%
\end{equation}
when $\mu_{n}=0$. In this particular situation, the solution of the equation
(\ref{intensity}) is given by the following Feynman-Kac path integral
formulae
\begin{equation}
\gamma_{n}(f)=\gamma_{0}(1)~\mathbb{E}\left(  f(X_{n})~\prod_{0\leq p<n}%
G_{p}(X_{p})\right)  \label{fkmod}%
\end{equation}
where $X_{n}$ stands for a Markov chain taking values in the state spaces
$E_{n}$ with initial distribution $\eta_{0}=\gamma_{0}/\gamma_{0}(1)$ and
Markov transitions $M_{n}$ (see for instance section 1.4.4.in~\cite{fk}).
These measure-valued equations have been studied at length in~\cite{fk}.

\subsection{Mean field particle interpretation}

\label{meanfieldsecin} The transport formula presented in (\ref{flotPhimc})
provides a natural interpretation of the probability distributions $\eta_{n}$
as the laws of a process $\overline{X}_{n}$ whose elementary transitions
$\overline{X}_{n}\leadsto\overline{X}_{n+1}$ depends on the distribution
$\eta_{n}=\mbox{\rm Law}(\overline{X}_{n})$ as well as on the current mass
$\gamma_{n}(1)$. In contrast to the more traditional McKean type nonlinear
Markov chains presented in \cite{fk}, the dependency on the mass process
induces a dependency of the whole sequence of measures $\eta_{p}$, from the
origin $p=0$ up to the current time $p=n$.

From now on, we will always assume that the mappings
\[
\left(  m,\left(  x^{i}\right)  _{1\leq i\leq N}\right)  \in\left(
\mathbb{R}_{+}\times E_{n}^{N}\right)  \mapsto K_{n+1,\left(  m,\frac{1}%
{N}\sum_{i=1}^{N}\delta_{x^{i}}\right)  }\left(  x,A_{n+1}\right)
\]
are measurable w.r.t. the product sigma fields on $(\mathbb{R}_{+}\times
E_{n}^{N})$, for any $n\geq0$, $N\geq1$, and $1\leq i\leq N$, and any
measurable subset $A_{n+1}\subset E_{n+1}$. In this situation, the mean field
particle interpretation of (\ref{flotPhimc}) is an $E_{n}^{N}$-valued
sequence
\index{$\xi_{n}^{(N)}$}
$\xi_{n}^{(N)}=\left(  \xi_{n}^{(N,i)}\right)  _{1\leq i\leq N}$ defined as
\begin{equation}
\left\{
\begin{array}
[c]{rcl}%
\gamma_{n+1}^{N}(1) & = & \gamma_{n}^{N}(1)~\eta_{n}^{N}(G_{n})+\mu_{n+1}(1)\\
&  & \\
\mathbb{P}\left(  \xi_{n+1}^{(N)}\in dx~\left\vert ~\mathcal{F}_{n}%
^{(N)}\right.  \right)  & = & \prod_{i=1}^{N}~K_{n+1,\left(  \gamma_{n}%
^{N}(1),\eta_{n}^{N}\right)  }(\xi_{n}^{(N,i)},dx^{i})
\end{array}
\right.  \label{meanfieldeta}%
\end{equation}
with the pair of occupation measures $\left(  \gamma_{n}^{N},\eta_{n}%
^{N}\right)  $ defined below
\[
\eta_{n}^{N}:=\frac{1}{N}\sum_{i=1}^{N}\delta_{\xi_{n}^{(N,i)}}\quad
\mbox{\rm and}\quad\gamma_{n}^{N}(dx):=\gamma_{n}^{N}(1)~\eta_{n}^{N}(dx)
\]
In the above displayed formula, $\mathcal{F}_{n}^{(N)}$ stands for the
$\sigma$-field generated by the random sequence $(\xi_{p}^{(N)})_{0\leq p\leq
n}$, and $dx=dx^{1}\times\ldots\times dx^{N}$ stands for an infinitesimal
neighborhood of a point $x=(x^{1},\ldots,x^{N})\in E_{n}^{N}$. The initial
system $\xi_{0}^{(N)}$ consists of $N$ independent and identically distributed
random variables with common law $\eta_{0}$. As usual, to simplify the
presentation, we will suppress the parameter $N$ when there is no possible
confusion, so that we write $\xi_{n}$ and $\xi_{n}^{i}$ instead of $\xi
_{n}^{(N)}$ and $\xi_{n}^{(N,i)}$.

In the above discussion, we have implicitly assumed that the quantities
$\mu_{n}(1)$ are known and that it is easy to sample from the probability
distribution $\overline{\mu}_{n}(dx):={\mu_{n}(dx)}/{\mu_{n}(1)}$. In
practice, we often need to resort to an additional approximation scheme to
approximate $\mu_{n}(1)$ and $\overline{\mu}_{n}$. This situation is discussed
in section~\ref{case3sec}. This additional level of approximation has
essentially a minimal impact on the properties of the particle approximation
scheme which can be analyzed using the same tools.

\subsection{Notation}

For the convenience of the reader, we end this introduction with some notation
used in the present article. We denote by $\mathcal{M}(E)$ the set of measures
on some measurable state space $(E,\mathcal{E})$ and we recall that
$\mathcal{P}(E)$ is the set of probability measures. We also denote
$\mathcal{B}(E)$ the Banach space of all bounded and measurable functions $f$
equipped with the uniform norm $\Vert f\Vert$ and $\mbox{Osc}_{1}(E)$ the
convex set of $\mathcal{E}$-measurable functions $f$ with oscillations
$\mbox{osc}(f)\leq1$ where $\mbox{osc}(f)=\underset{(x,y)\in E^{2}}{\sup
}\left\vert f\left(  x\right)  -f\left(  y\right)  \right\vert $.

We let $\mu(f)=\int~\mu(dx)~f(x)$ be the Lebesgue integral of a function
$f\in\mathcal{B}(E)$ with respect to a measure $\mu\in\mathcal{M}(E)$. We
recall that a bounded integral kernel $M(x,dy)$ from a measurable space
$(E,\mathcal{E})$ into an auxiliary measurable space $(E^{\prime}%
,\mathcal{E}^{\prime})$ is an operator $f\mapsto M(f)$ from $\mathcal{B}%
(E^{\prime})$ into $\mathcal{B}(E)$ such that the functions $x\mapsto
M(f)(x):=\int_{E^{\prime}}M(x,dy)f(y)$ are $\mathcal{E}$-measurable and
bounded for any $f\in\mathcal{B}(E^{\prime})$. The kernel $M$ also generates a
dual operator $\mu\mapsto\mu M$ from $\mathcal{M}(E)$ into $\mathcal{M}%
(E^{\prime})$ defined by $(\mu M)(f):=\mu(M(f))$. A Markov kernel is a
positive and bounded integral operator $M$ with $M(1)\left(  x\right)  =1$ for
any $x\in E$. Given a pair of bounded integral operators $(M_{1},M_{2})$, we
let $(M_{1}M_{2})$ be the composition operator defined by $(M_{1}%
M_{2})(f)=M_{1}(M_{2}(f))$. For time-homogenous state spaces, we denote by
$M^{k}=M^{k-1}M=MM^{k-1}$ the $k$-th composition of a given bounded integral
operator $M$, with $k\geq0$, with the convention $M^{0}=Id$ the identity
operator. We also use the notation
\[
M\left(  \left[  f_{1}-M(f_{1})\right]  \left[  f_{2}-M(f_{2})\right]
\right)  (x):=M\left(  \left[  f_{1}-M(f_{1})(x)\right]  \left[  f_{2}%
-M(f_{2})(x)\right]  \right)  (x)
\]
for some bounded functions $f_{1},f_{2}$.

We also denote the total variation norm on $\mathcal{M}(E)$ by $\Vert\mu
\Vert_{\mathrm{tv}}=\sup_{f\in\mbox{\rm\tiny
Osc}_{1}(E)}|\mu(f)|$. When the bounded integral operator $M$ has a constant
mass, that is $M(1)\left(  x\right)  =M(1)\left(  y\right)  $ for any
$(x,y)\in E^{2}$, the operator $\mu\mapsto\mu M$ maps $\mathcal{M}(E)$ into
$\mathcal{M}(E^{\prime})$. In this situation, we let $\beta(M)$ be the
Dobrushin coefficient of a bounded integral operator $M$ defined by the
following formula
\[
\beta(M):=\sup{\ \{\mbox{\rm osc}(M(f))\;;\;\;f\in\mbox{\rm Osc}_{1}(E)\}}%
\]
Given a positive function $G$ on $E$, we let $\Psi_{G}~:~\eta\in
\mathcal{P}(E)\mapsto\Psi_{G}(\eta)\in\mathcal{P}(E)$ be the Boltzmann-Gibbs
transformation defined by
\[
\Psi_{G}(\eta)(dx):=\frac{1}{\eta(G)}~G(x)~\eta(dx)
\]
We recall that $\Psi_{G}(\eta)$ can be expressed in terms of a Markov
transport equation
\begin{equation}
\eta S_{\eta}=\Psi_{G}(\eta) \label{TG}%
\end{equation}
for some selection type transition $S_{\eta}(x,dy)$. For instance, for any
$\epsilon\geq0$ s.t. $G(x)>\epsilon$ for any $x$, we notice that%
\[
\Psi_{(G-\epsilon)}(\eta)=\frac{\eta(G)}{\eta(G)-\epsilon}\left(  \Psi
_{(G)}(\eta)-\frac{\epsilon\eta}{\eta(G)}\right)
\]
so we can take
\begin{equation}
S_{\eta}(x,dy):=\frac{\epsilon}{\eta(G)}~\delta_{x}(dy)+\left(  1-\frac
{\epsilon}{\eta(G)}\right)  ~\Psi_{(G-\epsilon)}(\eta)(dy) \label{SGn1}%
\end{equation}
For $\epsilon=0$, we have $S_{\eta}(x,dy)=\Psi_{G}(\eta)(dy)$. We can also
choose
\begin{equation}
S_{\eta}(x,dy):=\epsilon G(x)~\delta_{x}(dy)+\left(  1-\epsilon G(x)\right)
~\Psi_{G}(\eta)(dy) \label{SGn2}%
\end{equation}
for any $\epsilon\geq0$ that may depend on the current measure $\eta$, and
s.t. $\epsilon G(x)\leq1$. For instance, we can choose $1/\epsilon$ to be the
$\eta$-essential supremum of $G$.

\section{Statement of the main results\label{sec:statementresults}}

At the end of section \ref{intensityflows}, we have seen that the evolution
equation (\ref{intensity}) coincides with that of a Feynman-Kac model
(\ref{fkmod}) for $\mu_{n}=0$. In this specific situation, the distributions
$\gamma_{n}$ are simply given by the recursive equation
\begin{equation}
\gamma_{n}=\gamma_{n-1}Q_{n}~\Longrightarrow~\forall0\leq p\leq n\qquad
\gamma_{n}=\gamma_{p}Q_{p,n}\quad\mbox{\rm with}\quad Q_{p,n}=Q_{p+1}\ldots
Q_{n-1}Q_{n} \label{fkmodell}%
\end{equation}
For $p=n$, we use the convention $Q_{n,n}=Id$. In addition, the nonlinear
semigroup associated to this sequence of distributions is given by
\begin{equation}
\eta_{n}(f)=\Phi_{p,n}(\eta_{p})(f):={\eta_{p}Q_{p,n}(f)}/{\eta_{p}Q_{p,n}%
(1)}={\eta_{p}\left(  Q_{p,n}(1)P_{p,n}(f)\right)  }/{\eta_{p}Q_{p,n}(1)}
\label{Phipn}%
\end{equation}
with the Markov kernel $P_{p,n}(x_{p},dx_{n})=Q_{p,n}(x_{p},dx_{n}%
)/{Q_{p,n}(x_{p},E_{n})}$. The analysis of the mean field particle
interpretations of such models has been studied in~\cite{fk}. Various
properties including contraction inequalities, fluctuations, large deviations
and concentration properties have been developed for this class of models. In
this context, the fluctuations properties as well as $\mathbb{L}_{r}$-mean
error estimates, including uniform estimates w.r.t. the time parameter are
often expressed in terms of two central parameters:
\begin{equation}
q_{p,n}=\sup_{x,y}\frac{Q_{p,n}(1)(x)}{Q_{p,n}(1)(y)}\quad\mbox{and}\quad
\beta(P_{p,n})=\sup_{x,y\in E_{p}}\Vert P_{p,n}(x,\mbox{\LARGE .})-P_{p,n}%
(y,\mbox{\LARGE .})\Vert_{\mathrm{tv}} \label{defqpn}%
\end{equation}
with the pair of Feynman-Kac semigroups $(P_{p,n},Q_{p,n})$ introduced in
(\ref{fkmodell}) and (\ref{Phipn}).

We also consider the pair of parameters $(g_{-}(n),g_{+}(n))$ defined below
\[
g_{-}(n)=\inf_{0\leq p<n}\inf_{E_{p}}G_{p}\leq\sup_{0\leq p<n}\sup_{E_{p}%
}G_{p}=g_{+}(n)
\]
We also write $g_{-/+}(n)$ to refer to both parameters. The first main
objective of this article is to extend some of these properties to models
where $\mu_{n}$ is non necessarily null. We illustrate our estimates in three
typical scenarios
\begin{equation}
1)\quad G=g_{-/+}=1\qquad2)\quad g_{+}<1\quad\mbox{\rm and}\quad3)\quad
g_{-}>1 \label{threetypes}%
\end{equation}
arising in time homogeneous models
\begin{equation}
(E_{n},G_{n},M_{n},\mu_{n},g_{-}(n),g_{+}(n))=(E,G,M,\mu,g_{-},g_{+})
\label{cashomog}%
\end{equation}
These three scenarios correspond to the case where, \emph{independently from
the additional spontaneous births}, the existing targets die or survive and
spawn in such a way that either their number remains constant ($G=g_{-/+}=1$),
decreases ($g_{+}<1$) or increases ($g_{-}>1$).

Our first main result concerns three different types of long time behavior for
these three types of models. This result can basically be stated as follows.

\begin{theo}
\label{premierth} For time homogeneous models (\ref{cashomog}), the limiting
behavior of $(\gamma_{n}(1),\eta_{n})$ in the three scenarios
(\ref{threetypes}) is as follows:

\begin{enumerate}
\item When $G(x)=1$ for any $x\in E$, we have
\[
\gamma_{n}(1)=\gamma_{0}(1)+\mu(1)~n\quad\mbox{and}\quad\Vert\eta_{n}%
-\eta_{\infty}\Vert_{\mathrm{tv}}=O\left(  \frac{1}{n}\right)
\]
when $M$ is chosen so that
\begin{equation}
\sum_{n\geq0}\sup_{x\in E}\Vert M^{n}(x,\mbox{\LARGE .})-\eta_{\infty}%
\Vert_{\mathrm{tv}}<\infty\quad\mbox{for some invariant measure
$\eta_{\infty}=\eta_{\infty}M$.} \label{hypoM}%
\end{equation}

\item When $g_{+}<1$, there exists a constant $c<\infty$ such that
\[
\forall f\in\mathcal{B}(E),\qquad\left\vert \gamma_{n}(f)-\gamma_{\infty
}(f)\right\vert \vee\left\vert \eta_{n}(f)-\eta_{\infty}(f)\right\vert \leq
c~g_{+}^{n}~\Vert f\Vert
\]
with the limiting measures%
\begin{equation}
\gamma_{\infty}(f):=\sum_{n\geq0}\mu Q^{n}(f)~\mbox{and}~\eta_{\infty
}(f):=\gamma_{\infty}(f)/\gamma_{\infty}(1) \label{limiting2}%
\end{equation}

\item When $g_{-}>1$ and there exist $k\geq1$ and $\epsilon>0$ such that
$M^{k}(x,\mbox{\LARGE .})\geq\epsilon~M^{k}(y,\mbox{\LARGE .})$ for any
$x,y\in E$ then the mapping $\Phi=\Phi_{n-1,n}$ introduced in (\ref{Phipn})
has a unique fixed point $\eta_{\infty}=\Phi(\eta_{\infty})$ and
\[
\lim_{n\rightarrow\infty}\frac{1}{n}\log{\gamma_{n}(1)}=\log{\eta_{\infty}%
(G)}\quad\mbox{and}\quad\Vert\eta_{n}-\eta_{\infty}\Vert_{\mathrm{tv}}\leq
c~e^{-\lambda n}%
\]
for some finite constant $c<\infty$ and some $\lambda>0$.
\end{enumerate}
\end{theo}

A more precise statement and a detailed proof of the above theorem can be
found in section~\ref{secprooftheo1}.

Our second main result concerns the convergence of the mean field particle
approximations presented in (\ref{meanfieldeta}). We provide rather sharp non
asymptotic estimates including uniform convergence results w.r.t. the time
parameter. Our results can be basically stated as follows.

\begin{theo}
\label{theointro} For any $n\geq0$, and any $N\geq1$, we have $\gamma_{n}(1)$
and $\gamma^{N}_{n}(1) \in I_{n}$ with the compact interval $I_{n}$ defined
below
\begin{equation}
\label{encadre}I_{n}:=\left[  m_{-}(n), m_{+}(n) \right]  \quad
\mbox{where}\quad m_{-/+}(n):= \sum_{p=0}^{n}\mu_{p}(1)g_{-/+}(n)^{(n-p)}%
\end{equation}

In addition, for any $r\geq1$, $f\in\mbox{\rm Osc}_{1}(E_{n})$, and any
$N\geq1$, we have
\begin{equation}
\sqrt{N}~\mathbb{E}\left(  \left\vert \left[  \eta_{n}^{N}-\eta_{n}\right]
(f)\right\vert ^{r}\right)  ^{\frac{1}{r}}\leq a_{r}~b_{n}\quad
\mbox{with}\quad b_{n}\leq\sum_{p=0}^{n}b_{p,n} \label{mean-error}%
\end{equation}
where $a_{r}<\infty$ stands for a constant whose value only depends on the
parameter $r$ and $b_{p,n}$ is the collection of constants given by
\begin{equation}
b_{p,n}:=2~\left(  1\wedge m_{p,n}\right)  ~q_{p,n}\left[  q_{p,n}%
~\beta(P_{p,n})+\sum_{p<q\leq n}\frac{c_{q,n}}{\sum_{p<r\leq n}c_{r,n}}%
~\beta(P_{q,n})\right]  \label{mean-error-ct}%
\end{equation}
with the pair of parameters
\[
m_{p,n}=m_{+}(p){\Vert Q_{p,n}(1)\Vert}/{\sum_{p<q\leq n}c_{q,n}}%
\quad\mbox{and}\quad c_{p,n}:=\mu_{p}Q_{p,n}(1)
\]

Furthermore, the particle measures $\gamma_{n}^{N}$ are unbiased, and for the
three scenarios (\ref{threetypes}) with time homogenous models s.t.
$M^{k}(x,\mbox{\LARGE .})\geq\epsilon~M^{k}(y,\mbox{\LARGE .})$, for any
$x,y\in E$ and some pair of parameters $k\geq1$ and $\epsilon>0$, the constant
$b_{n}$ in (\ref{mean-error}) can be chosen so that $\sup_{n\geq0}b_{n}%
<\infty$; in addition, we have the non asymptotic variance estimates for some
$d<\infty$, any $n\geq1$ and for any $N>1$
\begin{equation}
\mathbb{E}\left(  \left[  \frac{\gamma_{n}^{N}(1)}{\gamma_{n}(1)}-1\right]
^{2}\right)  \leq d~\frac{n+1}{N-1}~\left(  1+\frac{d}{N-1}\right)  ^{n-1}
\label{varestima}%
\end{equation}

\end{theo}

The non asymptotic estimates stated in the above theorem extend the one
presented in~\cite{cerou2008,fk} for Feynman-Kac type models (\ref{fkmod})
where $\mu_{n}=0$. For such models, the $\mathbb{L}_{r}$-mean error estimates
(\ref{mean-error}) are satisfied with the collection of parameters
$b_{p,n}:=2q_{p,n}^{2}~\beta(P_{p,n})$, with $p\leq n$. The extra terms in
(\ref{mean-error-ct}) are intimately related to $\mu_{n}$ whose effects in the
semigroup stability depend on the nature of $G_{n}$. We refer to
theorem~\ref{premierth}, section~\ref{secprooftheo1} and
section~\ref{stabsect}, for a discussion on three different behaviors in the
three cases presented in (\ref{threetypes}).

A direct consequence of this theorem is that it implies the almost sure
convergence results:
\[
\lim_{N\rightarrow\infty}\eta_{n}^{N}(f)=\eta_{n}(f)\quad\mbox{\rm and}\quad
\lim_{N\rightarrow\infty}\gamma_{n}^{N}(f)=\gamma_{n}(f)
\]
for any bounded function $f\in\mathcal{B}(E_{n})$.

Our last main result is a functional central limit theorem. We let $W_{n}^{N}$
be the centered random fields defined by the following formula
\begin{equation}
\eta_{n}^{N}=\eta_{n-1}^{N}K_{n,(\gamma_{n}^{N}(1),\eta_{n-1}^{N})}+\frac
{1}{\sqrt{N}}~W_{n}^{N}\,. \label{defWNn}%
\end{equation}
We also consider the pair of random fields
\[
V_{n}^{\eta,N}:=\sqrt{N}[\eta_{n}^{N}-\eta_{n}]\quad\mbox{\rm and}\quad
V_{n}^{\gamma,N}:=\sqrt{N}[\gamma_{n}^{N}-\gamma_{n}]
\]
For $n=0$, we use the convention $W_{0}^{N}=V_{0}^{\eta,N}$.

\begin{theo}
\label{theointrotcl} The sequence of random fields $(W_{n}^{N})_{n\geq0}$
converges in law, as $N$ tends to infinity, to the sequence of $n$
independent, Gaussian and centered random fields $(W_{n})_{n\geq0}$ with a
covariance function given for any $f,g\in\mathcal{B}(E_{n})$ and $n\geq0$ by
\begin{equation}%
\begin{array}
[c]{l}%
\mathbb{E}(W_{n}(f)W_{n}(g))\\
\\
=\eta_{n-1}K_{n,(\gamma_{n-1}(1),\eta_{n-1})}\left(  [f-K_{n,(\gamma
_{n-1}(1),\eta_{n-1})}(f)][g-K_{n,(\gamma_{n-1}(1),\eta_{n-1})}(g)]\right)
)\,.
\end{array}
\label{corr1}%
\end{equation}
In addition, the pair of random fields $V_{n}^{\gamma,N}$ and $V_{n}^{\eta,N}$
converge in law as $N\rightarrow\infty$ to a pair of centered Gaussian fields
$V_{n}^{\gamma}$ and $V_{n}^{\eta}$ defined by
\[
V_{n}^{\gamma}(f):=\sum_{p=0}^{n}\gamma_{p}(1)~W_{p}(Q_{p,n}(f))\quad
\mbox{and}\quad V_{n}^{\eta}(f):=V_{n}^{\gamma}\left(  \frac{1}{\gamma_{n}%
(1)}(f-\eta_{n}(f))\right)
\]

\end{theo}

The details of the proof of theorem~\ref{theointro} and
theorem~\ref{theointrotcl} can be found in section~\ref{secasympt}. The proof
of the non-asymptotic variance estimate (\ref{varestima}) is given in
section~\ref{secasymptun} dedicated to the convergence of the unnormalized
particle measures $\gamma_{n}^{N}$. The $\mathbb{L}_{r}$-mean error estimates
(\ref{mean-error}) and the fluctuation theorem~\ref{theointrotcl} are proved
in section~\ref{normalizedmod}. Under additional regularity conditions, we
conjecture that it is possible to obtain uniform estimates for
theorem~\ref{theointrotcl} but have not established it here.

The rest of the article is organized as follows.

In section~\ref{branchingsectsg}, we analyze the semigroup properties of the
total mass process $\gamma_{n}(1)$ and the sequence of probability
distributions $\eta_{n}$. This section is mainly concerned with the proof of
theorem~\ref{premierth}. The long time behavior of the total mass process is
discussed in section~\ref{totmasssec}, while the asymptotic behavior of the
probability distributions is discussed in section~\ref{secprooftheo1}. In
section~\ref{stabsect}, we develop a series of Lipschitz type functional
inequalities for uniform estimates w.r.t. the time parameter for the particle
approximation. In section~\ref{meanfieldsec}, we present the McKean models
associated to the sequence $(\gamma_{n}(1),\eta_{n})$ and their mean field
particle interpretations. Section~\ref{secasympt} is concerned with the
convergence analysis of these particle approximations. In
section~\ref{secasymptun}, we discuss the convergence of the approximations of
$\gamma_{n}(1)$, including their unbiasedness property and the non asymptotic
variance estimates presented in (\ref{varestima}). The proof of the
$\mathbb{L}_{r}$-mean error estimates (\ref{mean-error}) is presented in
section~\ref{normalizedmod}. The proof of the functional central limit
theorem~\ref{theointrotcl} is a more or less direct consequence of the
decomposition formulae presented in section~\ref{secasympt} and is just
sketched at the end of this very section.

\section{Semigroup analysis}

\label{branchingsectsg} The purpose of this section is to analyze the
semigroup properties of the intensity measure recursion (\ref{intensity}). We
establish a framework for the analysis of the long time behavior of these
measures and their particle approximations (\ref{meanfieldeta}). First, we
briefly recall some estimate of the quantities $(q_{p,n},\beta(P_{p,n}))$ in
terms of the potential functions $G_{n}$ and the Markov transitions $M_{n} $.
Further details on this subject can be found in~\cite{fk}, and in references therein.

We assume here that the following condition is satisfied for some $k\geq1$,
some collection of numbers $\epsilon_{p}\in(0,1)$
\begin{equation}
\hskip-1.5cm(M)_{k}\quad M_{p,p+k}(x_{p},\mbox{\LARGE .})\geq\epsilon
_{p}~M_{p,p+k}(y_{p},\mbox{\LARGE .})~~\mbox{\rm with}~~M_{p,p+k}%
=M_{p+1}M_{p+2}\ldots M_{p+k} \label{condMm}%
\end{equation}
for any time parameter $p$ and any pair of states $(x_{p},y_{p})\in E_{p}^{2}
$. It is well known that the mixing type condition $(M)_{k}$ is satisfied for
any aperiodic and irreducible Markov chains on finite spaces, as well as for
bi-Laplace exponential transitions associated with a bounded drift function
and for Gaussian transitions with a mean drift function that is constant
outside some compact domain. We introduce the following quantities
\begin{equation}
\delta_{p,n}:=\sup_{{}}{\prod_{p\leq q<n}\left(  {G_{q}(x_{q})}/{G_{q}(y_{q}%
)}\right)  }\quad\mbox{and}\quad\delta_{p}^{(k)}:=\delta_{p+1,p+k}
\label{fkhom2}%
\end{equation}
where the supremum is taken over all admissible pair of paths with transitions
$M_{q}$ where an admissible path $\left(  x_{p-1},x_{p+1},...,x_{n-1}\right)
$ is such that ${\prod_{p\leq q<n}M}_{q}\left(  x_{q-1},dx_{q}\right)  >0$.
Under the above conditions, we have~\cite[p. 140]{fk}
\begin{equation}
\beta(P_{p,p+n})\leq\prod_{l=0}^{\lfloor n/k\rfloor-1}\left(  1-\epsilon
_{p+lk}^{2}/\delta_{p+lk}^{(k)}\right)  \quad\mbox{and}\quad q_{p,p+n}%
\leq\delta_{p,p+k}/\epsilon_{p} \label{estenrs}%
\end{equation}
For time-homogeneous Feynman-Kac models we set $\epsilon:=\epsilon_{k}$ and
$\delta_{k}:=\delta_{0,k}$, for any $k\geq0$. Using this notation, the above
estimates reduce to~\cite[p. 142]{fk}
\begin{equation}
q_{p,p+n}\leq~\delta_{k}/\epsilon\quad\mbox{\rm
and}\quad\beta(P_{p,p+n})\leq\left(  1-\epsilon^{2}/\delta_{k-1}\right)
^{\lfloor n/k\rfloor} \label{estenr}%
\end{equation}

\subsection{Description of the models}

\label{totmasssec} The next proposition gives a Markov transport formulation
of $\Gamma_{n}$ introduced in (\ref{flotPhi}).

\begin{prop}
\label{proprefsg} For any $n\geq0$, we have the recursive formula
\begin{equation}
\label{phieq2}\left\{
\begin{array}
[c]{rcl}%
\gamma_{n+1}(1) & = & \gamma_{n}(1)~\eta_{n}(G_{n})+\mu_{n+1}(1)\\
&  & \\
\eta_{n+1} & = & \Psi_{G_{n}}(\eta_{n})M_{n+1,(\gamma_{n}(1),\eta_{n})}%
\end{array}
\right.
\end{equation}
with the collection of Markov transitions $M_{n+1,(m,\eta)}$ indexed by the
parameters $m\in\mathbb{R}_{+}$ and the probability measures $\eta
\in\mathcal{P }(E_{n})$ given below
\begin{equation}
\label{defM}M_{n+1,(m,\eta)}(x,dy):=\alpha_{n}\left(  m,\eta\right)
M_{n+1}(x,dy)+\left(  1-\alpha_{n}\left(  m,\eta\right)  \right)
~\overline{\mu}_{n+1}(dy)
\end{equation}
with the collection of $[0,1]$-parameters $\alpha_{n}\left(  m,\eta\right)  $
defined below
\[
\alpha_{n}\left(  m,\eta\right)  = \frac{m\eta(G_{n})}{m\eta(G_{n})+\mu
_{n+1}(1)}%
\]

\end{prop}

\noindent\mbox{\bf Proof:}\newline Observe that for any function
$f\in\mathcal{B}(E_{n+1})$, we have that
\[
\eta_{n+1}(f)=\frac{\gamma_{n}(G_{n}M_{n+1}(f))+\mu_{n+1}(f)}{\gamma_{n}%
(G_{n})+\mu_{n+1}(1)}=\frac{\gamma_{n}(1)~\eta_{n}(G_{n}M_{n+1}(f))+\mu
_{n+1}(f)}{\gamma_{n}(1)~\eta_{n}(G_{n})+\mu_{n+1}(1)}%
\]
from which we find that
\[
\eta_{n+1}=\alpha_{n}\left(  \gamma_{n}(1),\eta_{n}\right)  ~\Phi_{n+1}%
(\eta_{n})+\left(  1-\alpha_{n}\left(  \gamma_{n}(1),\eta_{n}\right)  \right)
~\overline{\mu}_{n+1}%
\]
From these observations, we prove (\ref{phieq2}). This ends the proof of the
proposition. \hfill\hbox{\vrule height 5pt width 5pt depth 0pt}\medskip
\newline

We let $\Gamma_{n+1}$ be the mapping from $\mathbb{R}_{+}\times\mathcal{P}%
(E_{n})$ into $\mathbb{R}_{+}\times\mathcal{P}(E_{n+1})$ given by
\begin{equation}
\Gamma_{n+1}(m,\eta)=\left(  \Gamma_{n+1}^{1}(m,\eta),\Gamma_{n+1}^{2}%
(m,\eta)\right)  \label{defGamma}%
\end{equation}
with the pair of transformations:
\[
\Gamma_{n+1}^{1}(m,\eta)=m~\eta(G_{n})+\mu_{n+1}(1)\quad\mbox{and}\quad
\Gamma_{n+1}^{2}(m,\eta)=\Psi_{G_{n}}(\eta)M_{n+1,(m,\eta)}%
\]
We also denote by $\left(  \Gamma_{p,n}\right)  _{0\leq p\leq n}$ the
corresponding semigroup defined by
\[
\forall0\leq p\leq n\qquad\Gamma_{p,n}=\Gamma_{p+1,n}\Gamma_{p+1}=\Gamma
_{n}\Gamma_{n-1}\ldots\Gamma_{p+1}%
\]
with the convention $\Gamma_{n,n}=Id$.

The following lemma collects some important properties of the sequence of
intensity measures $\gamma_{n}$.

\begin{lem}
\label{semigrouplin} For any $0\leq p\leq n$, we have the semigroup
decomposition
\begin{equation}
\label{gammaexplic}\gamma_{n}=\gamma_{p}Q_{p,n}+\sum_{p<q\leq n}\mu_{q}Q_{q,n}
\quad\mbox{and}\quad\gamma_{n}=\sum_{0\leq p\leq n}\mu_{p}Q_{p,n}%
\end{equation}
In addition, we also have the following formula
\begin{equation}
\label{prodform}\gamma_{n}(1)=\sum_{p=0}^{n}~\mu_{p}(1)~\prod_{p\leq q<n}%
~\eta_{q}(G_{q})
\end{equation}

\end{lem}

\noindent\mbox{\bf Proof:}\newline

The first pair of formulae are easily proved using a simple induction, and
recalling that $\gamma_{0}=\mu_{0}$. To prove the last assertion, we use an
induction on the parameter $n\geq0$. The result is obvious for $n=0$. We also
have by (\ref{intensity})
\[
\gamma_{n+1}(1)=\gamma_{n}Q_{n+1}(1)+\mu_{n+1}(1)=\gamma_{n}(G_{n})+\mu
_{n+1}(1)
\]
This implies
\begin{align*}
\gamma_{n+1}(1)  &  =\gamma_{n}(1)~\eta_{n}(G_{n})+\mu_{n+1}(1)\\
&  =\gamma_{n-1}(1)~\eta_{n-1}(G_{n-1})~\eta_{n}(G_{n})+\mu_{n}(1)~~\eta
_{n}(G_{n})+\mu_{n+1}(1)\\
&  =\ldots\\
&  =\gamma_{0}(1)~\prod_{p=0}^{n}~\eta_{p}(G_{p})~+\sum_{p=1}^{n+1}~\mu
_{p}(1)~\prod_{p\leq q\leq n}~\eta_{q}(G_{q})
\end{align*}
Recalling that $\gamma_{0}(dx_{0})=\mu_{0}(dx_{0})$, we prove (\ref{prodform}%
). This ends the proof of the lemma. \hfill
\hbox{\vrule height 5pt width 5pt depth
0pt}\medskip\newline

Using lemma~\ref{semigrouplin}, one proves that the semigroup $\Gamma_{p,n}$
satisfies the pair of formulae described below

\begin{prop}
\label{propsg} For any $0\leq p\leq n$, we have
\begin{align}
\Gamma_{p,n}^{1}(m,\eta)  &  =m~\eta Q_{p,n}(1)+\sum_{p<q\leq n}\mu_{q}%
Q_{q,n}(1)\label{sg1dec}\\
\Gamma_{p,n}^{2}(m,\eta)  &  = \alpha_{p,n}\left(  m,\eta\right)  ~ \Phi
_{p,n}(\eta)+\left(  1- \alpha_{p,n}\left(  m,\eta\right)  \right)
\sum_{p<q\leq n}\frac{c_{q,n}}{\sum_{p<r\leq n}c_{r,n}}~\Phi_{q,n}%
(\overline{\mu}_{q})\nonumber\\
&  \label{sg2dec}%
\end{align}
with the collection of parameters $c_{p,n}:=\mu_{p}Q_{p,n}(1)$ and the
$[0,1]$-valued parameters $\alpha_{p,n}\left(  m,\eta\right)  $ defined below
\begin{equation}
\label{eqalph}\alpha_{p,n}\left(  m,\eta\right)  =\frac{m\eta Q_{p,n}%
(1)}{m\eta Q_{p,n}(1)+\sum_{p<q\leq n}c_{q,n}} \leq\alpha^{\star}_{p,n}(m):=
1\wedge\left[  m\left\|  \frac{Q_{p,n}(1)}{\sum_{p<q\leq n}c_{q,n}}\right\|
\right]
\end{equation}

\end{prop}

One central question in the theory of spatial branching point processes is the
long time behavior of the total mass process $\gamma_{n}(1)$. Notice that
$\gamma_{n}(1)=\mathbb{E}(\mathcal{X}_{n}(1))$ is the expected size of the
$n$-th generation. For time homogeneous models with null spontaneous branching
$\mu_{n}=\mu=0$, the exponential growth of these quantities are related to the
logarithmic Lyapunov exponents of the semigroup $Q_{p,n}$. The prototype of
these models is the Galton-Watson branching process. In this context three
typical situations may occur: 1) $\gamma_{n}(1)$ remains constant and equals
to the initial mean number of individuals. 2) $\gamma_{n}(1)$ goes
exponentially fast to 0, 3) $\gamma_{n}(1)$ grows exponentially fast to infinity,

The analysis of spatial branching point processes with $\mu_{n}=\mu\neq0$
considered here is more involved. Loosely speaking, in the first situation
discussed above the total mass process is generally strictly increasing; while
in the second situation the additional mass injected in the system stabilizes
the total mass process. Before giving further details, by
lemma~\ref{semigrouplin} we observe $\gamma_{n}(1)\in I_{n}$, for any $n\geq
0$, with the compact interval $I_{n}$ defined in~\ref{encadre}.

We end this section with a more precise analysis of the effect of $\mu$ in the
three scenarios (\ref{threetypes}).

In the further developments of this section, we illustrate the stability
properties of the sequence of probability distributions $\eta_{n}$ in these
three scenarios.

\begin{enumerate}
\item When $G(x)=1$ for any $x\in E$, the total mass process $\gamma_{n}(1)$
grows linearly w.r.t. the time parameter and we have%
\begin{equation}
\gamma_{n}(1)=m_{-}(n)=m_{+}(n)=\gamma_{0}(1)+\mu(1)~n \label{caseG1}%
\end{equation}
Note that the estimates in (\ref{eqalph}) take the following form
\[
\alpha_{p,n}\left(  \gamma_{p}(1),\eta_{p}\right)  \leq\alpha_{p,n}^{\star
}(\gamma_{p}(1)):=1\wedge\frac{\gamma_{0}(1)+\mu(1)~p}{\mu(1)~(n-p)}%
\rightarrow_{(n-p)\rightarrow\infty}0
\]

\item When $g_{+}<1$, the total mass process $\gamma_{n}(1)$ is uniformly
bounded w.r.t. the time parameter. More precisely, we have that
\[
m_{-/+}(n)=g_{-/+}^{n}~\gamma_{0}(1)+\left(  1-g_{-/+}^{n}\right)  ~\frac
{\mu(1)}{1-g_{-/+}}%
\]
This yields the rather crude estimates
\begin{equation}
\gamma_{0}(1)\wedge\frac{\mu(1)}{1-g_{-}}\leq\gamma_{n}(1)\leq\gamma
_{0}(1)\vee\frac{\mu(1)}{1-g_{+}} \label{caseGl1}%
\end{equation}

We end this discussion with an estimate of the parameter $\alpha_{p,n}(m)$
given in (\ref{eqalph}). When the mixing condition $(M)_{k}$ stated in
(\ref{condMm}) is satisfied for some $k$ and some fixed parameters
$\epsilon_{p}=\epsilon$, using (\ref{estenr}) we prove that
\[
\sum_{p<r\leq n}\frac{\mu Q_{r,n}(1)}{Q_{p,r}(Q_{r,n}(1))}\geq\frac
{\epsilon\mu(1)}{\delta_{k}}~\sum_{p<r\leq n}\frac{1}{Q_{p,r}(1)}\geq
\frac{\epsilon\mu(1)}{\delta_{k}}~\frac{g_{+}^{-(n-p)}-1}{1-g_{+}}%
\]
from which we conclude that for any $n>p$ and any $m\in I_{p}$
\begin{align}
\alpha_{p,n}^{\star}(m)  &  \leq1\wedge\left[  m~g_{+}^{(n-p)}~\frac
{\delta_{k}~(1-g_{+})}{\epsilon\mu(1)(1-g_{+}^{(n-p)})}\right] \nonumber\\
&  \leq1\wedge\left[  m~g_{+}^{(n-p)}~{\delta_{k}}/{(\epsilon\mu(1))}\right]
\nonumber\\
&  \leq1\wedge\left[  \left(  \gamma_{0}(1)\vee\frac{\mu(1)}{1-g_{+}}\right)
~g_{+}^{(n-p)}~{\delta_{k}}/{(\epsilon\mu(1))}\right]  \rightarrow
_{(n-p)\rightarrow\infty}0 \label{refespi}%
\end{align}

\item When $g_{-}>1$, the total mass process $\gamma_{n}(1)$ grows
exponentially fast w.r.t. the time parameter and we can easily show that
\begin{equation}
g_{-}>1\Longrightarrow\gamma_{n}(1)\geq m_{-}(n)=\gamma_{0}(1)~g_{-}^{n}%
+\mu(1)~\frac{g_{-}^{n}-1}{g_{-}-1} \label{caseGb1}%
\end{equation}

\end{enumerate}

\subsection{Asymptotic properties}

\label{secprooftheo1} This section is concerned with the long time behavior of
the semigroups $\Gamma_{p,n}$ in the three scenarios discussed in
(\ref{caseG1}), (\ref{caseGl1}), and (\ref{caseGb1}). Our results are
summarized in theorem~\ref{premierth}. We consider time-homogeneous models
$(E_{n},G_{n},M_{n},\mu_{n})=(E,G,M,\mu)$.

\begin{enumerate}
\item When $G(x)=1$ for any $x\in E$, we have seen in (\ref{caseG1}) that
$\gamma_{n}(1)=\gamma_{0}(1)+\mu(1)~n$. In this particular situation, the
time-inhomogeneous Markov transitions $M_{n,(\gamma_{n-1}(1),\eta_{n-1}%
)}:=\overline{M}_{n}$ introduced in (\ref{phieq2}) are given by
\[
\overline{M}_{n}(x,dy)=\left(  1-\frac{\mu(1)}{\gamma_{0}(1)+n\mu(1)}\right)
~M(x,dy)+\frac{\mu(1)}{\gamma_{0}(1)+n\mu(1)}~\overline{\mu}(dy)
\]
This shows that $\eta_{n}=\mbox{\rm Law}(\overline{X}_{n})$ can be interpreted
as the distribution of the states $\overline{X}_{n}$ of a time inhomogeneous
Markov chain with transitions $\overline{M}_{n}$ and initial distribution
$\eta_{0}$. If we choose in (\ref{flotPhimc}) $K_{n+1,(\gamma_{n}(1),\eta
_{n})}=\overline{M}_{n+1}$, the $N$-particle model (\ref{meanfieldeta})
reduces to a series of $N$ independent copies of $\overline{X}_{n}$. In this
situation, the mapping $\Gamma_{0,n}^{2}$ is given by
\[
\Gamma_{0,n}^{2}(\gamma_{0}(1),\eta_{0}):=\frac{\gamma_{0}(1)}{\gamma
_{0}(1)+n\mu(1)}~~\eta_{0}M^{n}+\frac{n\mu(1)}{\gamma_{0}(1)+n\mu(1)}%
~~~\frac{1}{n}\sum_{0\leq p<n}\overline{\mu}M^{p}%
\]
The above formula shows that for a large time horizon $n$, the normalized
distribution flow $\eta_{n}$ is almost equal to $\frac{1}{n}\sum_{0\leq
p<n}\overline{\mu}M^{p}$. Let us assume that the Markov kernel $M$ is chosen
so that (\ref{hypoM}) is satisfied for some invariant measure $\eta_{\infty
}=\eta_{\infty}M$. In this case, for any starting measure $\gamma_{0}$, we
have
\[
\Vert\eta_{n}-\eta_{\infty}\Vert_{\mathrm{tv}}\leq\frac{\gamma_{0}(1)}%
{\gamma_{0}(1)+n\mu(1)}~\tau_{n}+\frac{n\mu(1)}{\gamma_{0}(1)+n\mu(1)}%
~~~\frac{1}{n}\sum_{0\leq p<n}\tau_{p}=O\left(  \frac{1}{n}\right)
\]
with $\tau_{n}=\sup_{x\in E}\Vert M^{n}(x,\mbox{\LARGE .})-\eta_{\infty}%
\Vert_{\mathrm{tv}}$. For instance, suppose the mixing condition $(M)_{k}$
presented in (\ref{condMm}) is met for some $k\geq1$ and $\epsilon>0$. In this
case, the above upper bound is satisfied with $\tau_{n}=(1-\epsilon)^{\lfloor
n/k\rfloor}$.

\item Consider the case where $g_{+}<1$. In this situation, the pair of
measures (\ref{limiting2}) are well defined. Furthermore, for any
$f\in\mathcal{B}(E)$ with $\Vert f\Vert\leq1$, we have the estimates
\begin{align*}
\left\vert \gamma_{n}(f)-\gamma_{\infty}(f)\right\vert  &  \leq\gamma
_{0}(1)~\eta_{0}Q^{n}(1)+\sum_{p\geq n}\mu Q^{p}(1)\\
&  \leq g_{+}^{n}~\left[  \gamma_{0}(1)+\mu(1)/(1-g_{+})\right]
\longrightarrow_{n\rightarrow\infty}0
\end{align*}
In addition, using the fact that $\gamma_{n}(1)\geq\mu(1)$, we find that for
any $f\in\mbox{\rm Osc}_{1}(E)$
\begin{align*}
\left\vert \eta_{n}(f)-\eta_{\infty}(f)\right\vert  &  \leq\frac{1}{\gamma
_{n}(1)}\left\vert \gamma_{n}[f-\eta_{\infty}(f)]-\gamma_{\infty}%
[f-\eta_{\infty}(f)]\right\vert \\
&  \leq g_{+}^{n}~\left[  \gamma_{0}(1)/\mu(1)+1/(1-g_{+})\right]
\longrightarrow_{n\rightarrow\infty}0
\end{align*}

\item Consider the case where $g_{-}>1$. We further assume that the mixing
condition $(M)_{k}$ presented in (\ref{condMm}) is met for some $k\geq1$ and
some fixed parameters $\epsilon_{p}=$ $\epsilon>0$. In this situation, it is
well known that the mapping $\Phi=\Phi_{n-1,n}$ introduced in (\ref{Phipn})
has a unique fixed point $\eta_{\infty}=\Phi(\eta_{\infty})$, and for any
initial distribution $\eta_{0}$, we have
\begin{equation}
\Vert\Phi_{0,n}(\eta_{0})-\eta_{\infty}\Vert_{\mathrm{tv}}\leq a~e^{-\lambda
~n} \label{sgcontract}%
\end{equation}
with
\[
\lambda=-\frac{1}{k}\log{\left(  1-\epsilon^{2}/\delta_{0,k-1}\right)  }%
\quad\mbox{\rm and}\quad a=1/\left(  1-\epsilon^{2}/\delta_{0,k-1}\right)
\]
as well as
\begin{equation}
\sup_{\eta\in\mathcal{P}(E)}\left\vert \frac{1}{n}\log{\eta Q^{n}(1)}%
-\log{\eta_{\infty}(G)}\right\vert \leq b/n \label{massasymp}%
\end{equation}
for some finite constant $b<\infty$. For a more thorough discussion on the
stability properties of the semigroup $\Phi_{0,n}$ and the limiting measures
$\eta_{\infty}$, we refer the reader to~\cite{fk}. Our next objective is to
transfer these stability properties to the one of the sequence $\eta_{n}$.
First, using (\ref{massasymp}), we readily prove that
\[
\lim_{n\rightarrow\infty}\frac{1}{n}\log{\gamma_{n}(1)}=\log{\eta_{\infty}(G)}%
\]
Next, we simplify the notation and we set $\alpha_{n}:=\alpha_{0,n}\left(
\gamma_{0}(1),\eta_{0}\right)  $ and $c_{n}:=c_{0,n}$. Using (\ref{sg2dec}),
we find that for any $n>1$
\[
a^{-1}~\Vert\eta_{n}-\eta_{\infty}\Vert_{\mathrm{tv}}\leq\alpha_{n}%
~e^{-\lambda n}+\left(  1-\alpha_{n}\right)  ~\sum_{0\leq p<n}\frac{c_{p}%
}{\sum_{0\leq q<n}c_{q}}~e^{-\lambda p}%
\]
Recalling that
\[
\mu(1)~g_{-}^{p}\leq c_{p}=\mu Q^{p}(1)\leq\mu(1)~g_{+}^{p}%
\]
we also obtain that
\begin{align}
\sum_{0\leq p<n}\frac{c_{p}}{\sum_{1\leq q<n}c_{q}}~e^{-\lambda p}  &
\leq\frac{1}{\left[  \sum_{0\leq q<n}c_{q}\right]  ^{1/r}}~\left[  \sum_{0\leq
p<n}c_{p}e^{-\lambda pr}\right]  ^{1/r}\nonumber\\
&  \leq\frac{1}{\left[  \sum_{0\leq q<n}g_{-}^{q}\right]  ^{1/r}}~\left[
\sum_{0\leq p<n}(e^{-\lambda r}g_{+})^{p}\right]  ^{1/r} \label{refcasgp1}%
\end{align}
for any $r\geq1$. We conclude that
\[
r>\frac{1}{\lambda}~\log{g_{+}}\Longrightarrow\sum_{0\leq p<n}\frac{c_{p}%
}{\sum_{0\leq q<n}c_{q}}~e^{-\lambda p}\leq g_{-}^{-(n-1)/r}/(1-e^{-\lambda
r}g_{+})^{1/r}%
\]
and therefore
\[
a^{-1}~\Vert\eta_{n}-\eta_{\infty}\Vert_{\mathrm{tv}}\leq e^{-\lambda n}%
+g_{-}^{-(n-1)/r}~/(1-e^{-\lambda r}g_{+})^{1/r}\rightarrow_{n\rightarrow
\infty}0
\]

\end{enumerate}

\subsection{Stability and Lipschitz regularity properties}

\label{stabsect} We describe in this section a framework that allows to
transfer the regularity properties of the Feynman-Kac semigroups $\Phi_{p,n} $
introduced in (\ref{Phipn}) to the ones of the semigroup $\Gamma_{p,n}$ of the
sequence $(\gamma_{n}(1),\eta_{n})$. Before proceeding we recall a lemma that
provides some weak Lipschitz type inequalities for the Feynman-Kac semigroup
$\Phi_{p,n}$ in terms of the Dobrushin contraction coefficient associated with
the Markov transitions $P_{p,n}$ introduced in (\ref{Phipn}). The details of
the proof of this result can be found in~\cite{fk} or in~\cite{ddj} (see Lemma
4.4. in~\cite{ddj}, or proposition 4.3.7 on page 146 in~\cite{fk}).

\begin{lem}
[\cite{ddj}]\label{lemFKPhi} For any $0\leq p\leq n$, any $\eta,\mu
\in\mathcal{P}(E_{p})$ and any $f\in\mbox{\rm Osc}_{1}(E_{n})$, we have
\begin{equation}
\left\vert \left[  \Phi_{p,n}(\mu)-\Phi_{p,n}(\eta)\right]  (f)\right\vert
\leq2~q_{p,n}^{2}~\beta(P_{p,n})~\left\vert (\mu-\eta)\mathcal{D}_{p,n,\eta
}(f)\right\vert \label{firstlip}%
\end{equation}
for a collection of functions $\mathcal{D}_{p,n,\eta}(f)\in\mbox{\rm Osc}_{1}%
(E_{p})$ whose values only depend on the parameters $(p,n,\eta)$.
\end{lem}

\begin{prop}
\label{regupropc2} For any $0\leq p\leq n$, any $\eta,\eta^{\prime}%
\in\mathcal{P}(E_{p})$ and any $f\in\mbox{\rm Osc}_{1}(E_{n})$, there exits a
collection of functions $\mathcal{D}_{p,n,\eta^{\prime}}(f)\in\mbox{\rm
Osc}_{1}(E_{p})$ whose values only depend on the parameters $(p,n,\eta)$ and
such that, for any $m\in I_{p}$, we have
\begin{equation}%
\begin{array}
[c]{l}%
\left\vert \left[  \Gamma_{p,n}^{2}(m,\eta)-\Gamma_{p,n}^{2}(m,\eta^{\prime
})\right]  (f)\right\vert \\
\\
\leq2~\alpha_{p,n}^{\star}~q_{p,n}\left[  q_{p,n}~\beta(P_{p,n})~\left\vert
(\eta-\eta^{\prime})\mathcal{D}_{p,n,\eta^{\prime}}(f)\right\vert +\beta
_{p,n}\left\vert (\eta-\eta^{\prime})h_{p,n,\eta^{\prime}}\right\vert \right]
\end{array}
\label{lipineq}%
\end{equation}
with the collection of functions $h_{p,n,\eta^{\prime}}=\frac{1}{2q_{p,n}%
}~\frac{Q_{p,n}(1)}{\eta^{\prime}Q_{p,n}(1)}\in\mbox{\rm Osc}_{1}(E_{p})$ and
the sequence of parameters $\epsilon_{p,n}$ and $\beta_{p,n}$ defined below
\begin{equation}
\alpha_{p,n}^{\star}:=\alpha_{p,n}^{\star}(m_{+}(p))\quad\mbox{and}\quad
\beta_{p,n}:=\sum_{p<q\leq n}\frac{c_{q,n}}{\sum_{p<r\leq n}c_{r,n}}%
~\beta(P_{q,n}) \label{defvarpi}%
\end{equation}

\end{prop}

Before getting into the details of the proof of proposition~\ref{regupropc2},
we illustrate some consequences of these weak functional inequalities for
time-homogeneous models $(E_{n},G_{n},M_{n},\mu_{n})=(E,G,M,\mu)$ in the three
scenarios discussed in (\ref{caseG1}), (\ref{caseGl1}), and (\ref{caseGb1}).

\begin{enumerate}
\item When $G(x)=1$ for any $x\in E$, we have
\[
\Phi_{p,n}(\eta)=\eta M^{(n-p)},\quad h_{p,n,\eta^{\prime}}=1/2\quad
c_{p,n}=\mu(1)\quad q_{p,n}=1\quad\alpha_{p,n}^{\star}\leq1
\]
Let us assume that there exist $a<\infty$ and $0<\lambda<\infty$ such that
$\beta(M^{n})\leq ae^{-\lambda n}$ for any $n\geq0$. In this situation, we
prove using (\ref{lipineq}) that
\[%
\begin{array}
[c]{l}%
\left\vert \left[  \Gamma_{p,n}^{2}(m,\eta)-\Gamma_{p,n}^{2}(m,\eta^{\prime
})\right]  (f)\right\vert \leq2ae^{-\lambda(n-p)}~\left\vert (\mu
-\eta)\mathcal{D}_{p,n,\eta^{\prime}}(f)\right\vert
\end{array}
\]

\item When $g_{+}<1$ and when the mixing condition $(M)_{k}$ stated in
(\ref{condMm}) is satisfied for some $k$ and some fixed parameters
$\epsilon_{p}=\epsilon$, we have seen in (\ref{refespi}) that
\[
\sup_{m\in I_{p}}\alpha_{p,n}^{\star}(m)\leq1\wedge\left(  d~g_{+}%
^{(n-p)}\right)  \quad\mbox{\rm
with}\quad d=\left(  (\gamma_{0}(1)/\mu(1))\vee(1-g_{+})^{-1}\right)
{\delta_{0,k}}{\epsilon^{-1}}%
\]
Furthermore, using the estimates given in (\ref{estenrs}) and (\ref{estenr}),
we also have that
\[
q_{p,n}\leq\delta_{k}/\epsilon\qquad\beta_{p,n}\leq1\quad\mbox{\rm
and}\quad\beta(P_{p,n})\leq a~e^{-\lambda~(n-p)}\quad\mbox{\rm
with $(a,\lambda)$ given in (\ref{sgcontract})}
\]
In this situation, we prove using (\ref{lipineq}) that
\[%
\begin{array}
[c]{l}%
\left\vert \left[  \Gamma_{p,n}^{2}(m,\eta)-\Gamma_{p,n}^{2}(m,\eta^{\prime
})\right]  (f)\right\vert \\
\\
\leq2~\left[  1\wedge\left(  d~g_{+}^{(n-p)}\right)  \right]  ~(\delta
_{k}/\epsilon)\left[  (\delta_{k}/\epsilon)~a~e^{-\lambda(n-p)}~\left\vert
(\mu-\eta)\mathcal{D}_{p,n,\eta^{\prime}}(f)\right\vert +\left\vert (\mu
-\eta)h_{p,n,\eta^{\prime}}\right\vert \right]
\end{array}
\]
Notice that for $(n-p)\geq\log{(d)}/\log{(1/g_{+})}$, this yields
\[%
\begin{array}
[c]{l}%
\left\vert \left[  \Gamma_{p,n}^{2}(m,\eta)-\Gamma_{p,n}^{2}(m,\eta^{\prime
})\right]  (f)\right\vert \\
\\
\leq a_{0}~e^{-\lambda_{0}(n-p)}~\left\vert (\mu-\eta)\mathcal{D}%
_{p,n,\eta^{\prime}}(f)\right\vert +a_{1}~e^{-\lambda_{1}(n-p)}\left\vert
(\mu-\eta)h_{p,n,\eta^{\prime}}\right\vert
\end{array}
\]
with
\[
a_{0}=2ad(\delta_{k}/\epsilon)^{2}\quad a_{1}=2d(\delta_{k}/\epsilon
)\quad\lambda_{0}=\lambda+\log{(1/g_{+})}\quad\mbox{\rm and}\quad\lambda
_{1}=\log{(1/g_{+})}%
\]

\item When $g_{-}>1$ and when the mixing condition $(M)_{k}$ presented in
(\ref{condMm}) is met for some $k$ and some fixed parameters $\epsilon
_{p}=\epsilon>0$, then we use the fact that
\[
\alpha_{p,n}^{\star}\leq1\quad q_{p,n}\leq\delta_{k}/\epsilon\quad\mbox{\rm
and}\quad\beta(P_{p,n})\leq a~e^{-\lambda(n-p)}\quad\mbox{\rm
with $(a,\lambda)$ given in (\ref{sgcontract})}
\]
Arguing as in (\ref{refcasgp1}), we prove that for any $r>\frac{1}{\lambda
}~\log{g_{+}}$
\[
\beta_{p,n}\leq g_{-}^{-(n-p-1)/r}/(1-e^{-\lambda r}g_{+})^{1/r}%
\]
from which we conclude that
\[%
\begin{array}
[c]{l}%
\left\vert \left[  \Gamma_{p,n}^{2}(m,\eta)-\Gamma_{p,n}^{2}(m,\eta^{\prime
})\right]  (f)\right\vert \\
\\
\leq a_{0}~e^{-\lambda_{0}(n-p)}~\left\vert (\mu-\eta)\mathcal{D}%
_{p,n,\eta^{\prime}}(f)\right\vert +a_{1}~e^{-\lambda_{1}(n-p)}\left\vert
(\mu-\eta)h_{p,n,\eta^{\prime}}\right\vert
\end{array}
\]
with
\[
a_{0}=2a(\delta_{k}/\epsilon)^{2}\quad a_{1}=2g_{-}^{r}(\delta_{k}%
/\epsilon)/(1-e^{-\lambda r}g_{+})^{1/r}\quad\lambda_{0}=\lambda
\quad\mbox{\rm and}\quad\lambda_{1}=\log{(g_{-})}%
\]

\end{enumerate}

Now, we come to the proof of proposition~\ref{regupropc2}.

\textbf{Proof of proposition~\ref{regupropc2}:}

First, we observe that
\[%
\begin{array}
[c]{l}%
\Gamma_{p,n}^{2}(m,\eta)-\Gamma_{p,n}^{2}(m^{\prime},\eta^{\prime})\\
\\
=\alpha_{p,n}\left(  m,\eta\right)  ~\left[  \Phi_{p,n}(\eta)-\sum_{p<q\leq
n}\frac{c_{q,n}}{\sum_{p<r\leq n}c_{r,n}}~\Phi_{q,n}(\overline{\mu}%
_{q})\right] \\
\\
\qquad\qquad-\alpha_{p,n}\left(  m^{\prime},\eta^{\prime}\right)  ~\left[
\Phi_{p,n}(\eta^{\prime})-\sum_{p<q\leq n}\frac{c_{q,n}}{\sum_{p<r\leq
n}c_{r,n}}~\Phi_{q,n}(\overline{\mu}_{q})\right]
\end{array}
\]
Using the following decomposition
\begin{equation}
ab-a^{\prime}b^{\prime}=a^{\prime}(b-b^{\prime})+(a-a^{\prime})b^{\prime
}+(a-a^{\prime})(b-b^{\prime}) \label{ab}%
\end{equation}
which is valid for any $a,a^{\prime},b,b^{\prime}\in\mathbb{R}$, we prove
that
\begin{equation}%
\begin{array}
[c]{l}%
\Gamma_{p,n}^{2}(m,\eta)-\Gamma_{p,n}^{2}(m^{\prime},\eta^{\prime})\\
\\
=\alpha_{p,n}\left(  m^{\prime},\eta^{\prime}\right)  \left[  \Phi_{p,n}%
(\eta)-\Phi_{p,n}(\eta^{\prime})\right] \\
\\
\qquad+\left[  \Phi_{p,n}(\eta^{\prime})-\sum_{p<q\leq n}\frac{c_{q,n}}%
{\sum_{p<r\leq n}c_{r,n}}~\Phi_{q,n}(\overline{\mu}_{q})\right]  \left[
\alpha_{p,n}\left(  m,\eta\right)  -\alpha_{p,n}\left(  m^{\prime}%
,\eta^{\prime}\right)  \right] \\
\\
\hskip5.3cm+\left[  \alpha_{p,n}\left(  m,\eta\right)  -\alpha_{p,n}\left(
m^{\prime},\eta^{\prime}\right)  \right]  \left[  \Phi_{p,n}(\eta)-\Phi
_{p,n}(\eta^{\prime})\right]
\end{array}
\label{gam2}%
\end{equation}

For $m=m^{\prime}$, using (\ref{gam2}) we find that
\[%
\begin{array}
[c]{l}%
\Gamma_{p,n}^{2}(m,\eta)-\Gamma_{p,n}^{2}(m,\eta^{\prime})\\
\\
=\alpha_{p,n}\left(  m,\eta\right)  ~\left[  \Phi_{p,n}(\eta)-\Phi_{p,n}%
(\eta^{\prime})\right] \\
\\
\qquad+\left[  \Phi_{p,n}(\eta^{\prime})-\sum_{p<q\leq n}\frac{c_{q,n}}%
{\sum_{p<r\leq n}c_{r,n}}~\Phi_{q,n}(\overline{\mu}_{q})\right]  \left[
\alpha_{p,n}\left(  m,\eta\right)  -\alpha_{p,n}\left(  m,\eta^{\prime
}\right)  \right]
\end{array}
\]
We also notice that
\[
\alpha_{p,n}\left(  m,\eta\right)  =\frac{1}{1+{\mu_{p,n}}/{\left[  m\eta
Q_{p,n}(1)\right]  }}%
\]
from which we easily prove that%
\[%
\begin{array}
[c]{l}%
\alpha_{p,n}\left(  m,\eta\right)  -\alpha_{p,n}\left(  m^{\prime}%
,\eta^{\prime}\right) \\
\\
=\frac{\mu_{p,n}}{\mu_{p,n}+m\eta Q_{p,n}(1)}\frac{1}{\mu_{p,n}+m^{\prime}%
\eta^{\prime}Q_{p,n}(1)}\left[  m\eta Q_{p,n}(1)-m^{\prime}\eta^{\prime
}Q_{p,n}(1)\right]
\end{array}
\]
and therefore
\[
\alpha_{p,n}\left(  m,\eta\right)  -\alpha_{p,n}\left(  m,\eta^{\prime
}\right)  \newline\newline=\left(  \alpha_{p,n}\left(  m,\eta^{\prime}\right)
(1-\alpha_{p,n}\left(  m,\eta\right)  )\right)  ~\left[  \eta-\eta^{\prime
}\right]  \left(  \frac{Q_{p,n}(1)}{\eta^{\prime}Q_{p,n}(1)}\right)
\]
The proof of $\alpha_{p,n}\left(  m,\eta\right)  \leq\alpha_{p,n}^{\star}(m)$
is elementary. From the above decomposition, we prove the following upper
bounds
\[
\left\vert \alpha_{p,n}\left(  m,\eta\right)  -\alpha_{p,n}\left(
m,\eta^{\prime}\right)  \right\vert \leq\alpha_{p,n}^{\star}(m)~\left\vert
\left[  \eta-\eta^{\prime}\right]  \left(  \frac{Q_{p,n}(1)}{\eta^{\prime
}Q_{p,n}(1)}\right)  \right\vert
\]
and
\[%
\begin{array}
[c]{l}%
\left\vert \left[  \Gamma_{p,n}^{2}(m,\eta)-\Gamma_{p,n}^{2}(m,\eta^{\prime
})\right]  (f)\right\vert \\
\\
\leq\alpha_{p,n}^{\star}(m)\left[  \left\vert \left[  \Phi_{p,n}(\eta
)-\Phi_{p,n}(\eta^{\prime})\right]  (f)\right\vert \right. \\
\\
\left.  \qquad+\left\vert \left[  \eta-\eta^{\prime}\right]  \left(
\frac{Q_{p,n}(1)}{\eta^{\prime}Q_{p,n}(1)}\right)  \right\vert ~\left\vert
\sum_{p<q\leq n}\frac{c_{q,n}}{\sum_{p<r\leq n}c_{r,n}}~\left[  \Phi
_{q,n}(\overline{\mu}_{q})-\Phi_{q,n}\left(  \Phi_{p,q}(\eta^{\prime})\right)
\right]  (f)\right\vert \right]
\end{array}
\]
This yields
\[%
\begin{array}
[c]{l}%
\left\vert \left[  \Gamma_{p,n}^{2}(m,\eta)-\Gamma_{p,n}^{2}(m,\eta^{\prime
})\right]  (f)\right\vert \\
\\
\leq\alpha_{p,n}^{\star}(m)\left[  \left\vert \left[  \Phi_{p,n}(\eta
)-\Phi_{p,n}(\eta^{\prime})\right]  (f)\right\vert +\beta_{p,n}~\left\vert
\left[  \eta-\eta^{\prime}\right]  \left(  \frac{Q_{p,n}(1)}{\eta^{\prime
}Q_{p,n}(1)}\right)  \right\vert \right]
\end{array}
\]
The last formula comes from the fact that
\[
\beta(P_{q,n}):=\sup_{\nu,\nu^{\prime}\in\mathcal{P}(E_{q})}\left\Vert
\Phi_{q,n}(\nu)-\Phi_{q,n}(\nu^{\prime})\right\Vert _{\mathrm{tv}}%
\]
The proof of this result can be found in~\cite{fk} (proposition 4.3.1 on page
134). The end of the proof is now a direct consequence of lemma~\ref{lemFKPhi}%
. This ends the proof of the proposition. \hfill\hbox{\vrule
height 5pt width 5pt depth 0pt}\medskip\newline

\section{Mean field particle approximations}

\label{meanfieldsec}

\subsection{McKean particle interpretations}

\label{nonlinsecc} In proposition~\ref{proprefsg}, the evolution equation
(\ref{phieq2}) of the sequence of probability measures $\eta_{n}\leadsto
\eta_{n+1}$ is a combination of an updating type transition $\eta_{n}%
\leadsto\Psi_{G_{n}}(\eta_{n})$ and an integral transformation w.r.t. a Markov
transition $M_{n+1,(\gamma_{n}(1),\eta_{n})}$ that depends on the current
total mass $\gamma_{n}(1)$ and the current probability distribution $\eta_{n}%
$. The operator $M_{n+1,(\gamma_{n}(1),\eta_{n})}$ defined in (\ref{defM}) is
a mixture of the Markov transition $M_{n+1}$ and the spontaneous birth
normalized measure $\overline{\mu}_{n+1}$. We let $S_{n,\eta_{n}}$ be any
Markov transition from $E_{n}$ into itself satisfying
\[
\Psi_{G_{n}}(\eta_{n})=\eta_{n}S_{n,\eta_{n}}%
\]
The choice of these transitions is not unique. We can choose for instance one
of the collection of transitions presented in (\ref{TG}), (\ref{SGn1}) and
(\ref{SGn2}). Further examples of McKean acceptance-rejection type transitions
can also be found in section 2.5.3 in~\cite{fk}. By construction, we have the
recursive formula
\begin{equation}
\eta_{n+1}=\eta_{n}K_{n+1,(\gamma_{n}(1),\eta_{n})}\quad\mbox{\rm with}\quad
K_{n+1,(\gamma_{n}(1),\eta_{n})}=S_{n,\eta_{n}}M_{n+1,(\gamma_{n}(1),\eta
_{n})} \label{systnl}%
\end{equation}
with the auxiliary total mass evolution equation
\begin{equation}
\gamma_{n+1}(1)=\gamma_{n}(1)~\eta_{n}(G_{n})+\mu_{n+1}(1) \label{systmass}%
\end{equation}
As already mentioned in section \ref{secintro}, the sequence of probability
distributions $\eta_{n}$ can be interpreted as the distributions of the states
$\overline{X}_{n}$ of a process defined, conditional upon $\left(  \gamma
_{n}(1),\eta_{n}\right)  $, by the elementary transitions
\[
\mathbb{P}\left(  \overline{X}_{n+1}\in dx~|~\overline{X}_{n}\right)
=K_{n,(\gamma_{n}(1),\eta_{n})}\left(  \overline{X}_{n},dx\right)
\quad\mbox{\rm with}\quad\eta_{n}=\mbox{\rm Law}(\overline{X}_{n})
\]
Next, we define the mean field particle interpretations of the sequence
$(\gamma_{n}(1),\eta_{n})$ given in (\ref{systnl}) and (\ref{systmass}).
First, mimicking formula (\ref{systmass}) we set
\[
\gamma_{n+1}^{N}(1):=\gamma_{n}^{N}(1)~\eta_{n}^{N}(G_{n})+\mu_{n+1}%
(1)\quad\mbox{\rm and}\quad\gamma_{n}^{N}(f)=\gamma_{n}^{N}(1)~\times~\eta
_{n}^{N}(f)
\]
for any $f\in\mathcal{B}(E_{n})$, with the initial measure $\gamma_{0}%
^{N}=\gamma_{0}$. It is important to notice that
\[
\gamma_{n}^{N}(1)=\gamma_{0}(1)~\prod_{0\leq q<n}~\eta_{q}^{N}(G_{q}%
)+\sum_{p=1}^{n}~\mu_{p}(1)~\prod_{p\leq q<n}~\eta_{q}^{N}(G_{q}%
)\Longrightarrow\gamma_{n}^{N}(1)\in I_{n}%
\]
The mean field particle interpretation of the nonlinear measure valued model
(\ref{systnl}) is an $E_{n}^{N}$-valued process $\xi_{n}$ with elementary
transitions defined in (\ref{meanfieldeta}) and (\ref{systnl}). By
construction, the particle evolution is a simple combination of a selection
and a mutation genetic type transition
\[
\xi_{n}\leadsto\widehat{\xi}_{n}=(\widehat{\xi}_{n}^{i})_{1\leq i\leq
N}\leadsto\xi_{n+1}%
\]
During the selection transitions $\xi_{n}\leadsto\widehat{\xi}_{n}$, each
particle $\xi_{n}^{i}\leadsto\widehat{\xi}_{n}^{i}$ evolves according to the
selection type transition $S_{n,\eta_{n}^{N}}(\xi_{n}^{i},dx)$. During the
mutation stage, each of the selected particles $\widehat{\xi}_{n}^{i}%
\leadsto\xi_{n+1}^{i}$ evolves according to the transition
\[
M_{n+1,(\gamma_{n}^{N}(1),\eta_{n}^{N})}(x,dy):=\alpha_{n}\left(  \gamma
_{n}^{N}(1),\eta_{n}^{N}\right)  M_{n+1}(x,dy)+\left(  1-\alpha_{n}\left(
\gamma_{n}^{N}(1),\eta_{n}^{N}\right)  \right)  ~\overline{\mu}_{n+1}(dy)
\]

\subsection{Asymptotic behavior}

\label{secasympt} This section is mainly concerned with the proof of
theorem~\ref{theointro}. In section~\ref{secasymptun}, we discuss the
unibiasedness property of the particle measures $\gamma_{n}^{N}$ and their
convergence properties towards $\gamma_{n}$, as the number of particles $N$
tends to infinity. We mention that the proof of the non asymptotic variance
estimates (\ref{varestima}) is simpler than the one provided in a recent
article by the second author with F. Cérou and A. Guyader~\cite{cerou2008}.
Section~\ref{normalizedmod} is concerned with the convergence and the
fluctuations of the occupation measures $\eta_{n}^{N}$ around their limiting
measures $\eta_{n}$.

\subsubsection{Intensity measures}

\label{secasymptun} We start this section with a simple unbiasedness property.
Recall that $\mathcal{F}_{p}^{(N)}$ stands for the $\sigma$-field generated by
the random sequence $(\xi_{k}^{(N)})_{0\leq k\leq p}$.

\begin{prop}
For any $0\leq p\leq n$, and any $f\in\mathcal{B}(E_{n})$, we have
\begin{equation}
\mathbb{E}\left(  \gamma_{n+1}^{N}(f)\left\vert ~\mathcal{F}_{p}^{(N)}\right.
\right)  =\gamma_{p}^{N}Q_{p,n+1}(f)+\sum_{p<q\leq n+1}\mu_{q}Q_{q,n+1}(f)
\label{unbiasp}%
\end{equation}
In particular, we have the unbiasedness property: $\mathbb{E}\left(
\gamma_{n}^{N}(f)\right)  =\gamma_{n}(f)$.
\end{prop}

\noindent\mbox{\bf Proof:}\newline By construction of the particle model, for
any $f\in\mathcal{B}(E_{n})$ we have
\[
\mathbb{E}\left(  \eta_{n+1}^{N}(f)\left\vert ~\mathcal{F}_{n}^{(N)}\right.
\right)  =\eta_{n}^{N}K_{n+1,(\gamma_{n}^{N}(1),\eta_{n}^{N})}(f)=\Gamma
_{n+1}^{2}\left(  \gamma_{n}^{N}(1),\eta_{n}^{N}\right)  (f)
\]
with the second component $\Gamma_{n+1}^{2}$ of the transformation
$\Gamma_{n+1}$ introduced in~\ref{defGamma}. Using the fact that
\[
\Gamma_{n+1}^{2}\left(  \gamma_{n}^{N}(1),\eta_{n}^{N}\right)  (f)=\frac
{\gamma_{n}^{N}(1)~\eta_{n}^{N}(Q_{n+1}(f))+\mu_{n+1}(f)}{\gamma_{n}%
^{N}(1)~\eta_{n}^{N}(Q_{n+1}(1))+\mu_{n+1}(1)}=\frac{\gamma_{n}^{N}%
(Q_{n+1}(f))+\mu_{n+1}(f)}{\gamma_{n}^{N}(Q_{n+1}(1))+\mu_{n+1}(1)}%
\]
and
\[
\gamma_{n+1}^{N}(1)=\gamma_{n}^{N}(1)~\eta_{n}^{N}(G_{n})+\mu_{n+1}%
(1)=\gamma_{n}^{N}(Q_{n+1}(1))+\mu_{n+1}(1)
\]
we prove that
\begin{align*}
\mathbb{E}\left(  \gamma_{n+1}^{N}(f)\left\vert ~\mathcal{F}_{n}^{(N)}\right.
\right)   &  =\mathbb{E}\left(  \gamma_{n+1}^{N}(1)~\eta_{n+1}^{N}%
(f)\left\vert ~\mathcal{F}_{n}^{(N)}\right.  \right)  =\gamma_{n+1}%
^{N}(1)~\mathbb{E}\left(  \eta_{n+1}^{N}(f)\left\vert ~\mathcal{F}_{n}%
^{(N)}\right.  \right) \\
&  =\gamma_{n}^{N}(Q_{n+1}(f))+\mu_{n+1}(f)
\end{align*}
This also implies that
\begin{align*}
\mathbb{E}\left(  \gamma_{n+1}^{N}(f)\left\vert ~\mathcal{F}_{n-1}%
^{(N)}\right.  \right)   &  =\mathbb{E}\left(  \gamma_{n}^{N}(Q_{n+1}%
(f))\left\vert ~\mathcal{F}_{n-1}^{(N)}\right.  \right)  +\mu_{n+1}(f)\\
&  =\gamma_{n-1}^{N}(Q_{n}Q_{n+1}(f))+\mu_{n}(Q_{n+1}(f))+\mu_{n+1}(f)
\end{align*}
Iterating the argument one proves (\ref{unbiasp}). The end of the proof is now
clear. \hfill\hbox{\vrule height 5pt width 5pt depth 0pt}\medskip\newline

The next theorem provides a key martingale decomposition and a rather crude
non asymptotic variance estimate.

\begin{theo}
\label{letheo1} For any $n\geq0$ and any function $f\in\mathcal{B}(E_{n})$, we
have the decomposition
\begin{equation}
\sqrt{N}~\left[  \gamma_{n}^{N}-\gamma_{n}\right]  (f)=\sum_{p=0}^{n}%
\gamma_{p}^{N}(1)~W_{p}^{N}(Q_{p,n}(f)) \label{martingaledec}%
\end{equation}
In addition, if the mixing condition $(M)_{k}$ presented in (\ref{condMm}) is
met for some $k\geq1$ and some constant parameters $\epsilon_{p}=\epsilon>0$,
then we have for any $N>1$ and any $n\geq1$
\begin{equation}
\mathbb{E}\left(  \left[  \frac{\gamma_{n}^{N}(1)}{\gamma_{n}(1)}-1\right]
^{2}\right)  \leq\frac{n+1}{N-1}~\frac{\delta_{k}^{2}}{\epsilon^{2}}\left(
1+\frac{\delta_{k}^{2}}{\epsilon^{2}(N-1)}\right)  ^{n-1} \label{varestim}%
\end{equation}

\end{theo}

Before presenting the proof of this theorem, we would like to make a couple of
comments. On the one hand, we observe that the unbiasedness property follows
directly from the decomposition (\ref{martingaledec}). On the other hand,
using Kintchine's inequality, for any $r\geq1$, $p\geq1$, and any
$f\in\mbox{\rm Osc}_{1}(E_{n})$ we have the almost sure estimates
\[
\sqrt{N}~\mathbb{E}\left(  \left\vert W_{p}^{N}(f)\right\vert ^{r}~\left\vert
\mathcal{F}_{p-1}^{(N)}\right.  \right)  ^{\frac{1}{r}}\leq a_{r}%
\]

A detailed proof of these estimates can be found in~\cite{fk}, see also lemma
7.2 in~\cite{bdd} for a simpler proof by induction on the parameter $N$. From
this elementary observation, and recalling that $\gamma_{n}^{N}(1)\in I_{n}$
for any $n\geq0$, we find that
\[
\sqrt{N}~\mathbb{E}\left(  \left\vert \left[  \gamma_{n}^{N}-\gamma
_{n}\right]  (f)\right\vert ^{r}\right)  ^{\frac{1}{r}}\leq a_{r}~b_{n}%
\]
for some finite constant $b_{n}$ whose values only depend on the time
parameter $n$.

Now, we present the proof of theorem~\ref{letheo1}.

\textbf{Proof of theorem~\ref{letheo1}:}

We use the decomposition:
\[
\gamma_{n+1}^{N}(f)-\gamma_{n+1}(f)=\left[  \gamma_{n+1}^{N}(f)-\mathbb{E}%
\left(  \gamma_{n+1}^{N}(f)\left\vert ~\mathcal{F}_{n}^{(N)}\right.  \right)
\right]  +\left[  \mathbb{E}\left(  \gamma_{n+1}^{N}(f)\left\vert
~\mathcal{F}_{n}^{(N)}\right.  \right)  -\gamma_{n+1}(f)\right]
\]
By (\ref{unbiasp}), we find that
\[
\gamma_{n+1}^{N}(f)-\mathbb{E}\left(  \gamma_{n+1}^{N}(f)\left\vert
~\mathcal{F}_{n}^{(N)}\right.  \right)  =\gamma_{n+1}^{N}(f)-\left[
\gamma_{n}^{N}(Q_{n+1}(f))+\mu_{n+1}(f)\right]
\]
Since we have
\begin{align*}
\gamma_{n}^{N}(Q_{n+1}(1))+\mu_{n+1}(1)  &  =\gamma_{n}^{N}(G_{n})+\mu
_{n+1}(1)\\
&  =\gamma_{n}^{N}(1)~\eta_{n}^{N}(G_{n})+\mu_{n+1}(1)=\gamma_{n+1}^{N}(1)
\end{align*}
this implies that
\begin{align*}
\gamma_{n+1}^{N}(f)-\left[  \gamma_{n}^{N}(Q_{n+1}(f))+\mu_{n+1}(f)\right]
&  =\gamma_{n+1}^{N}(1)~\left[  \eta_{n+1}^{N}(f)-\frac{\left[  \gamma_{n}%
^{N}(Q_{n+1}(f))+\mu_{n+1}(f)\right]  }{\left[  \gamma_{n}^{N}(Q_{n+1}%
(1))+\mu_{n+1}(1)\right]  }\right] \\
&  =\gamma_{n+1}^{N}(1)~\left[  \eta_{n+1}^{N}(f)-\eta_{n}^{N}K_{n+1,(\gamma
_{n}^{N}(1),\eta_{n}^{N})}(f)\right]
\end{align*}
and therefore
\[
\gamma_{n+1}^{N}(f)-\mathbb{E}\left(  \gamma_{n+1}^{N}(f)\left\vert
~\mathcal{F}_{n}^{(N)}\right.  \right)  =\gamma_{n+1}^{N}(1)~\left[
\eta_{n+1}^{N}(f)-\eta_{n}^{N}K_{n+1,(\gamma_{n}^{N}(1),\eta_{n}^{N}%
)}(f)\right]
\]
Finally, we observe that
\[
\mathbb{E}\left(  \gamma_{n+1}^{N}(f)\left\vert ~\mathcal{F}_{n}^{(N)}\right.
\right)  -\gamma_{n+1}(f)=\gamma_{n}^{N}(Q_{n+1}(f))-\gamma_{n}(Q_{n+1}(f))
\]
from which we find the recursive formula
\[
\left[  \gamma_{n+1}^{N}-\gamma_{n+1}\right]  (f)=\gamma_{n+1}^{N}(1)~\left[
\eta_{n+1}^{N}-\eta_{n}^{N}K_{n+1,(\gamma_{n}^{N}(1),\eta_{n}^{N})}\right]
(f)+\left[  \gamma_{n}^{N}-\gamma_{n}\right]  (Q_{n+1}(f))
\]
The end of the proof of (\ref{martingaledec}) is now obtained by a simple
induction on the parameter $n$.

Now, we come to the proof of (\ref{varestim}). Using the fact that
\begin{align*}
\mathbb{E}\left(  \gamma_{p}^{N}(1)W_{p}^{N}(f^{(1)})~\gamma_{q}^{N}%
(1)W_{q}^{N}(f^{(2)})\right)   &  =\mathbb{E}\left(  \gamma_{p}^{N}%
(1)\gamma_{q}^{N}(1)W_{p}^{N}(f^{(1)})~\mathbb{E}\left(  W_{q}^{N}%
(f^{(2)})~|~\mathcal{F}_{q-1}^{N}\right)  \right) \\
&  =0
\end{align*}
for any $0\leq p<q\leq n$, and any $f^{(1)}\in\mathcal{B}(E_{p})$, and
$f^{(2)}\in\mathcal{B}(E_{q})$, we prove that
\[
N~\mathbb{E}\left(  \left[  \gamma_{n}^{N}(1)-\gamma_{n}(1)\right]
^{2}\right)  =\sum_{p=0}^{n}\mathbb{E}\left(  \gamma_{p}^{N}(1)^{2}%
~\mathbb{E}\left(  W_{p}^{N}(Q_{p,n}(1))^{2}|\mathcal{F}_{p-1}^{N}\right)
\right)
\]
Notice that
\begin{equation}
\frac{1}{\gamma_{n}(1)^{2}}=\frac{1}{\gamma_{p}(1)^{2}}~\frac{1}{\eta
_{p}(Q_{p,n}(1))^{2}}~\left(  \frac{\gamma_{p}(Q_{p,n}(1))}{\gamma_{n}%
(1)}\right)  ^{2}\leq\alpha_{p,n}^{\star}(\gamma_{p}(1))^{2}\frac{1}%
{\gamma_{p}(1)^{2}}~\frac{1}{\eta_{p}(Q_{p,n}(1))^{2}} \label{eqreff}%
\end{equation}
The r.h.s. estimate comes from the fact that
\[
\frac{\gamma_{p}(Q_{p,n}(1))}{\gamma_{n}(1)}=\frac{\gamma_{p}(1)~\eta
_{p}(Q_{p,n}(1))}{\gamma_{p}(1)~\eta_{p}(Q_{p,n}(1))+\sum_{p<q\leq n}\mu
_{q}Q_{q,n}(1)}=\alpha_{p,n}\left(  \gamma_{p}(1),\eta_{p}\right)  \leq
\alpha_{p,n}^{\star}(\gamma_{p}(1))
\]
Using the above decompositions, we readily prove that
\[
N~\mathbb{E}\left(  \left[  \frac{\gamma_{n}^{N}(1)}{\gamma_{n}(1)}-1\right]
^{2}\right)  \leq\sum_{p=0}^{n}\alpha_{p,n}^{\star}(\gamma_{p}(1))^{2}%
~\mathbb{E}\left(  \left(  \frac{\gamma_{p}^{N}(1)}{\gamma_{p}(1)}\right)
^{2}~\mathbb{E}\left(  W_{p}^{N}(\overline{Q}_{p,n}(1))^{2}|\mathcal{F}%
_{p-1}^{N}\right)  \right)
\]
with
\[
\overline{Q}_{p,n}(1)=\overline{Q}_{p,n}(1)/\eta_{p}(Q_{p,n}(1))\leq q_{p,n}%
\]
We set
\[
U_{n}^{N}:=\mathbb{E}\left(  \left[  \frac{\gamma_{n}^{N}(1)}{\gamma_{n}%
(1)}-1\right]  ^{2}\right)  \quad\mbox{\rm
then we find that}\quad N~U_{n}^{N}\leq a_{n}+\sum_{p=0}^{n}b_{p,n}~U_{p}^{N}%
\]
with the parameters
\[
a_{n}:=\sum_{p=0}^{n}\left(  q_{p,n}\alpha_{p,n}^{\star}(\gamma_{p}(1)\right)
^{2}\quad\mbox{\rm and}\quad b_{p,n}:=\left(  q_{p,n}\alpha_{p,n}^{\star
}(\gamma_{p}(1)\right)  ^{2}%
\]
Using the fact that $b_{n,n}\leq1$, we prove the following recursive equation
\[
U_{n}^{N}\leq a_{n}^{N}+\sum_{0\leq p<n}b_{p,n}^{N}~U_{p}^{N}\quad\mbox{\rm
with}\quad a_{n}^{N}:=\frac{a_{n}}{N-1}\quad\mbox{\rm and}\quad b_{p,n}%
^{N}:=\frac{b_{p,n}}{N-1}%
\]
Using an elementary proof by induction on the time horizon $n$, we prove the
following inequality:
\[
U_{n}^{N}\leq\left[  \sum_{p=1}^{n}a_{p}^{N}\sum_{e\in\langle p,n\rangle}%
b^{N}(e)\right]  +\left[  \sum_{e\in\langle0,n\rangle}b^{N}(e)\right]
~U_{0}^{N}%
\]
In the above display, $\langle p,n\rangle$ stands for the set of all integer
valued paths $e=\left(  e(l)\right)  _{0\leq l\leq k}$ of a given length $k$
from $p$ to $n$
\[
e_{0}=p<e_{1}<\ldots<e_{k-1}<e_{k}=n\quad\mbox{\rm and}\quad b^{N}%
(e)=\prod_{1\leq l\leq k}b_{e(l-1),e(l)}^{N}%
\]
We have also used the convention $b^{N}(\emptyset)=\prod_{\emptyset}=1$ and
$\langle n,n\rangle=\{\emptyset\}$, for $p=n$. Recalling that $\gamma_{0}%
^{N}=\gamma_{0}$, we conclude that
\[
U_{n}^{N}\leq\sum_{p=1}^{n}a_{p}^{N}\sum_{e\in\langle p,n\rangle}b^{N}(e)
\]
We further assume that the mixing condition $(M)_{k}$ presented in
(\ref{condMm}) is met for some parameters $k\geq1$, and some constant
parameters $\epsilon_{p}=\epsilon>0$. In this case, we use the fact that
\[
\alpha_{p,n}^{\star}\leq1\quad\mbox{\rm and}\quad q_{p,n}\leq\delta
_{k}/\epsilon
\]
to prove that
\[
\sup_{0\leq p\leq n}a_{p}^{N}\leq(n+1)~(\delta_{k}/\epsilon)^{2}%
/(N-1)\quad\mbox{\rm and}\quad\sup_{0\leq p\leq n}b_{p,n}^{N}\leq(\delta
_{k}/\epsilon)^{2}/(N-1)
\]
Using these rather crude estimates, we find that
\[
U_{n}^{N}\leq a_{n}^{N}+\sum_{0<p<n}a_{p}^{N}\sum_{l=1}^{(n-p)}\left(
\begin{array}
[c]{c}%
n-p-1\\
l-1
\end{array}
\right)  \left(  \frac{\delta_{k}^{2}}{\epsilon^{2}(N-1)}\right)  ^{l}%
\]
and therefore
\begin{align*}
U_{n}^{N}  &  \leq\frac{(n+1)}{(N-1)}~\frac{\delta_{k}^{2}}{\epsilon^{2}%
}\left(  1+\frac{\delta_{k}^{2}}{\epsilon^{2}(N-1)}\sum_{0<p<n}\left(
1+\left(  \frac{\delta_{k}^{2}}{\epsilon^{2}(N-1)}\right)  \right)
^{n-p-1}\right) \\
&  =\frac{(n+1)}{(N-1)}~\frac{\delta_{k}^{2}}{\epsilon^{2}}\left(
1+\frac{\delta_{k}^{2}}{\epsilon^{2}(N-1)}\right)  ^{n-1}%
\end{align*}

This ends the proof of the theorem. \hfill\hbox{\vrule height 5pt width 5pt
depth 0pt}\medskip\newline

\subsubsection{Probability distributions}

\label{normalizedmod} This section is mainly concerned with the proof of the
$\mathbb{L}_{r}$-mean error estimates stated in (\ref{mean-error}). We use the
decomposition
\begin{align}
\left(  \gamma_{n}^{N}(1),\eta_{n}^{N}\right)  -\left(  \gamma_{n}(1),\eta
_{n}\right)   &  =\left[  \Gamma_{0,n}\left(  \gamma_{0}^{N}(1),\eta_{0}%
^{N}\right)  -\Gamma_{0,n}\left(  \gamma_{0}(1),\eta_{0}\right)  \right]
\nonumber\\
&  +\sum_{p=1}^{n}\left[  \Gamma_{p,n}\left(  \gamma_{p}^{N}(1),\eta_{p}%
^{N}\right)  -\Gamma_{p-1,n}\left(  \gamma_{p-1}^{N}(1),\eta_{p-1}^{N}\right)
\right]  \label{keydecomp}%
\end{align}
to prove that
\[%
\begin{array}
[c]{l}%
\eta_{n}^{N}-\eta_{n}\\
\\
=\left[  \Gamma_{0,n}^{2}\left(  \gamma_{0}^{N}(1),\eta_{0}^{N}\right)
-\Gamma_{0,n}^{2}\left(  \gamma_{0}(1),\eta_{0}\right)  \right]  +\sum
_{p=1}^{n}\left[  \Gamma_{p,n}^{2}\left(  \gamma_{p}^{N}(1),\eta_{p}%
^{N}\right)  -\Gamma_{p-1,n}^{2}\left(  \gamma_{p-1}^{N}(1),\eta_{p-1}%
^{N}\right)  \right]
\end{array}
\]
Using the fact that
\[
\Gamma_{p-1,n}(m,\eta)=\Gamma_{p,n}\left(  \Gamma_{p}(m,\eta)\right)
\Rightarrow\Gamma_{p-1,n}^{2}(m,\eta)=\Gamma_{p,n}^{2}\left(  \Gamma
_{p}(m,\eta)\right)
\]
we readily check that
\begin{align*}
\Gamma_{p}\left(  \gamma_{p-1}^{N}(1),\eta_{p-1}^{N}\right)   &  =\left(
\gamma_{p-1}^{N}(1)\eta_{p-1}^{N}(G_{p-1})+\mu_{p}(1),\Psi_{G_{p-1}}\left(
\eta_{p-1}^{N}\right)  M_{p,(\gamma_{p-1}^{N}(1),\eta_{p-1}^{N})}\right) \\
&  =\left(  \gamma_{p}^{N}(1),\eta_{p-1}^{N}K_{p,(\gamma_{p-1}^{N}%
(1),\eta_{p-1}^{N})}\right)
\end{align*}
Since we have $\gamma_{0}^{N}(1)=\mu_{0}(1)=\gamma_{0}(1)$, one concludes
that
\begin{align*}
\eta_{n}^{N}-\eta_{n}  &  =\left[  \Gamma_{0,n}^{2}\left(  \gamma_{0}%
(1),\eta_{0}^{N}\right)  -\Gamma_{0,n}^{2}\left(  \gamma_{0}(1),\eta
_{0}\right)  \right] \\
&  \qquad+\sum_{p=1}^{n}\left[  \Gamma_{p,n}^{2}\left(  \gamma_{p}^{N}%
(1),\eta_{p}^{N}\right)  -\Gamma_{p,n}^{2}\left(  \gamma_{p}^{N}(1),\eta
_{p-1}^{N}K_{p,(\gamma_{p-1}^{N}(1),\eta_{p-1}^{N})}\right)  \right]
\end{align*}
Using the fact that $\gamma_{p}^{N}(1)\in I_{p}$, for any $p\geq0$, the end of
the proof is a direct consequence of lemma~\ref{regupropc2} and Kintchine
inequality. The proof of the uniform convergence estimates stated in the end
of theorem~\ref{theointro} are a more or less direct consequence of the
functional inequalities derived at the end of section~\ref{stabsect}. The end
of the proof of the theorem~\ref{theointro} is now completed.

We end this section with the fluctuations properties of the $N$-particle
approximation measures $\gamma_{n}^{N}$ and $\eta_{n}^{N}$ around their
limiting values. Using the type of arguments as those used in the proof of the
functional central limit theorem, theorem 3.3 in \cite{dmrio}, we can prove
that the sequence $(W_{n}^{N})_{n\geq0}$ defined in (\ref{defWNn}) converges
in law, as $N$ tends to infinity, to the sequence of $n$ independent, Gaussian
and centered random fields $(W_{n})_{n\geq0}$ with a covariance function given
in (\ref{corr1}). Using the decompositions (\ref{martingaledec}) and
\[
\eta_{n}^{N}(f)-\eta_{n}(f)=\frac{\gamma_{n}(1)}{\gamma_{n}^{N}(1)}~\left(
[\gamma_{n}^{N}-\gamma_{n}]\left(  \frac{1}{\gamma_{n}(1)}(f-\eta
_{n}(f))\right)  \right)
\]
by the continuous mapping theorem, we deduce the functional central limit
theorem~\ref{theointrotcl}.

\section{Particle approximations of spontaneous birth measures}

\label{case3sec} Assume that the spontaneous birth measures $\mu_{n}$ are
chosen so that $\mu_{n}\ll\lambda_{n}$ for some reference probability measures
$\lambda_{n}$ and that the Radon Nikodim derivatives $H_{n}=d\mu_{n}%
/d\lambda_{n}$ are bounded. For any $n\geq0$, we let $\lambda_{n}^{N^{\prime}%
}:=\frac{1}{{N^{\prime}}}\sum_{i=1}^{N^{\prime}}\delta_{\zeta_{n}^{i}}$ be the
empirical measure associated with $N^{\prime}$ independent and identically
distributed random variables $\left(  \zeta_{n}^{i}\right)  _{1\leq i\leq N}$
with common distribution $\lambda_{n}$. We also denote by $\mu_{n}^{N^{\prime
}}$ the particle spontaneous birth measures defined below
\[
\forall n\geq0\qquad\mu_{n}^{N^{\prime}}(dx_{n}):=H_{n}(x_{n})~\lambda
_{n}^{N^{\prime}}(dx_{n})
\]
In this notation, the initial distribution $\eta_{0}$ and the initial mass
$\gamma_{0}$ are approximated by the weighted occupation measure $\eta
_{0}^{N^{\prime}}:=\Psi_{H_{0}}(\lambda_{0}^{N^{\prime}})$ and $\gamma
_{0}^{N^{\prime}}(1):=\lambda_{0}^{N^{\prime}}(H_{0})$.

We let $\widetilde{\gamma}_{n}^{N^{\prime}}$ and $\widetilde{\eta}%
_{n}^{N^{\prime}}$ the random measures defined as $\gamma_{n}$ and $\eta_{n}$
by replacing in (\ref{intensity}) the measures $\mu_{n}$ by the random
measures $\mu_{n}^{N^{\prime}}$, for any $n\geq0$; that is, we have that
\[
\widetilde{\gamma}_{n}^{N^{\prime}}=\widetilde{\gamma}_{n-1}^{N^{\prime}}%
Q_{n}+\mu_{n}^{N^{\prime}}\quad\mbox{\rm and}\quad\widetilde{\eta}%
_{n}^{N^{\prime}}(f_{n})=\widetilde{\gamma}_{n}^{N^{\prime}}(f_{n}%
)/\widetilde{\gamma}_{n}^{N^{\prime}}(1)
\]
for any $f_{n}\in\mathcal{B}(E_{n})$. By construction, using the same
arguments as the ones we used in the proof of (\ref{gammaexplic}), we have
\[
\widetilde{\gamma}_{n}^{N^{\prime}}=\sum_{0\leq p\leq n}\mu_{p}^{N^{\prime}%
}Q_{p,n}%
\]
This yields for any $f\in\mathcal{B}(E_{n})$ the decomposition
\[
\left[  \widetilde{\gamma}_{n}^{N^{\prime}}-\gamma_{n}\right]  (f)=\sum_{0\leq
p\leq n}\left[  \mu_{p}^{N^{\prime}}-\mu_{p}\right]  Q_{p,n}(f)=\sum_{0\leq
p\leq n}\left[  \lambda_{p}^{N^{\prime}}-\lambda_{p}\right]  \left(
H_{p}~Q_{p,n}(f)\right)
\]
Several estimates can be derived from these formulae, including $\mathbb{L}%
_{p}$-mean error bounds, functional central limit theorems, empirical process
convergence, as well as sharp exponential concentration inequalities. For
instance, we have the unbiasedness property
\[
\mathbb{E}\left(  \widetilde{\gamma}_{n}^{N^{\prime}}(f)\right)  =\gamma
_{n}(f)
\]
and the variance estimate
\[
N~\mathbb{E}\left(  \left[  \widetilde{\gamma}_{n}^{N^{\prime}}(f)-\gamma
_{n}(f)\right]  ^{2}\right)  =\sum_{0\leq p\leq n}\lambda_{p}\left[  \left(
H_{p}Q_{p,n}(f)-\lambda_{p}(H_{p}Q_{p,n}(f))\right]  ^{2}\right)
\]
Using the same arguments as the ones we used in (\ref{eqreff}), we prove the
following rather crude upper bound
\begin{align*}
N~\mathbb{E}\left(  \left[  \frac{\widetilde{\gamma}_{n}^{N^{\prime}}%
(f)}{\gamma_{n}(1)}-\eta_{n}(f)\right]  ^{2}\right)   &  \leq\sum_{0\leq p\leq
n}\alpha_{p,n}^{\star}(\gamma_{p}(1))^{2}\frac{1}{\gamma_{p}(1)^{2}}~\frac
{\mu_{p}\left(  H_{p}Q_{p,n}(f)^{2}\right)  }{\eta_{p}(Q_{p,n}(1))^{2}}\\
&  \leq\sum_{0\leq p\leq n}\alpha_{p,n}^{\star}(\gamma_{p}(1))^{2}\frac
{1}{\gamma_{p}(1)^{2}}~\Vert H_{p}\Vert~\mu_{p}(1)~q_{p,n}^{2}%
\end{align*}
We illustrate these variance estimates for time homogeneous models
$(E_{n},G_{n},H_{n},M_{n},\mu_{n})=(E,G,H,M,\mu)$, in the three situations
discussed in (\ref{caseG1}), (\ref{caseGl1}), and (\ref{caseGb1}). We further
assume that the mixing condition $(M)_{k}$ presented in (\ref{condMm}) is met
for some parameters $k\geq1$, and some $\epsilon>0$. In this case, we use the
fact that $q_{p,n}\leq\delta_{k}/\epsilon$, to prove that
\[
N~\mathbb{E}\left(  \left[  \frac{\widetilde{\gamma}_{n}^{N^{\prime}}%
(f)}{\gamma_{n}(1)}-\eta_{n}(f)\right]  ^{2}\right)  \leq c\sum_{0\leq p\leq
n}\left[  \alpha_{p,n}^{\star}(\gamma_{p}(1))/{\gamma_{p}(1)}\right]  ^{2}%
\]
with some constant $c:=\left(  \Vert H\Vert~\mu(1)~(\delta_{k}/\epsilon
)^{2}\right)  $.

\begin{enumerate}
\item When $G(x)=1$ for any $x\in E$, we have $\gamma_{p}(1)=\gamma_{0}%
(1)+\mu(1)~p$. Recalling that $\alpha_{p,n}^{\star}(\gamma_{p}(1))\leq1$, we
prove the uniform estimates
\[
N~\sup_{n\geq0}{\mathbb{E}\left(  \left[  \frac{\widetilde{\gamma}%
_{n}^{N^{\prime}}(f)}{\gamma_{n}(1)}-\eta_{n}(f)\right]  ^{2}\right)  }\leq
c\sum_{p\geq0}(\gamma_{0}(1)+\mu(1)~p)^{-2}%
\]

\item When $g_{+}<1$ and when the mixing condition $(M)_{k}$ stated in
(\ref{condMm}) is satisfied, we have seen in (\ref{refespi}) that
\[
\alpha_{p,n}^{\star}(\gamma_{p}(1))\leq1\wedge\left(  d_{1}~g_{+}%
^{(n-p)}\right)  \quad\mbox{\rm and}\quad\inf_{n}\gamma_{n}(1)\geq d_{2}%
\]
for some finite constants $d_{1}<\infty$ and $d_{2}>0$. From previous
calculations, we prove the following uniform variance estimates
\[
N~\sup_{n\geq0}{\mathbb{E}\left(  \left[  \frac{\widetilde{\gamma}%
_{n}^{N^{\prime}}(f)}{\gamma_{n}(1)}-\eta_{n}(f)\right]  ^{2}\right)  }%
\leq(c/d_{2}^{2})\sum_{p\geq0}\left[  1\wedge\left(  d_{1}^{2}~g_{+}%
^{2p}\right)  \right]
\]

\item When $g_{-}>1$ we have seen in (\ref{caseGb1}) that $\gamma_{n}(1)\geq
d~g_{-}^{n}$ for any $n\geq n_{0}$, for some finite constant $d<\infty$ and
some $n_{0}\geq1$ so
\[
N~\sup_{n\geq0}{\mathbb{E}\left(  \left[  \frac{\widetilde{\gamma}%
_{n}^{N^{\prime}}(f)}{\gamma_{n}(1)}-\eta_{n}(f)\right]  ^{2}\right)  }\leq
c\left(  \sum_{0\leq p\leq n_{0}}\gamma_{p}(1)^{-2}+d\sum_{n\geq n_{0}}%
g_{-}^{-2n}\right)
\]

\end{enumerate}

\end{document}